\documentclass[12pt, a4paper, reqno, oneside]{article}

\usepackage{amssymb,amsmath}
\DeclareUnicodeCharacter{0335}{}
\usepackage{tikz-cd}
\usepackage{authblk} 
\usepackage{physics} 
\usepackage{tikz}
\usepackage{titlesec} 
\usetikzlibrary{babel}
\usepackage{adjustbox}
\usepackage{url}
\usepackage{xcolor}
\usepackage{yhmath}
\usepackage{eso-pic}
\usepackage[T1]{fontenc}
\usepackage[utf8]{inputenc}
\usepackage[spanish,english]{babel}
\usepackage{csquotes}
\usepackage{pst-node}
\usepackage[english]{babel}
\usepackage{latexsym}
\usepackage{euscript}
\usepackage{mathtools}
\usepackage[all]{xy}
\usepackage{mathrsfs}
\usepackage[percent]{overpic}
\usepackage{setspace}
\usepackage{graphicx}
\usepackage{pdfsync}
\usepackage{amsthm}
\usepackage{thmtools}
\usepackage{yfonts}
\usetikzlibrary{calc}
\usetikzlibrary{decorations.pathmorphing}

\tikzset{curve/.style={settings={#1},to path={(\tikztostart)
    .. controls ($(\tikztostart)!\pv{pos}!(\tikztotarget)!\pv{height}!270:(\tikztotarget)$)
    and ($(\tikztostart)!1-\pv{pos}!(\tikztotarget)!\pv{height}!270:(\tikztotarget)$)
    .. (\tikztotarget)\tikztonodes}},
    settings/.code={\tikzset{quiver/.cd,#1}
        \def\pv##1{\pgfkeysvalueof{/tikz/quiver/##1}}},
    quiver/.cd,pos/.initial=0.35,height/.initial=0}

\tikzset{tail reversed/.code={\pgfsetarrowsstart{tikzcd to}}}
\tikzset{2tail/.code={\pgfsetarrowsstart{Implies[reversed]}}}
\tikzset{2tail reversed/.code={\pgfsetarrowsstart{Implies}}}
\tikzset{no body/.style={/tikz/dash pattern=on 0 off 1mm}}

\usepackage[colorlinks=true,linktocpage=true, pagebackref=false, citecolor=black,linkcolor=black]{hyperref}
\normalbaselines
\makeatletter

\def\ps@myfancy{\let\@mkboth\markboth
 \def\@evenhead{\vbox{\hsize\textwidth 
 \hbox to \textwidth{\sf\mdseries\thepage 
 \rule[-.6ex]{0mm}{2mm} \hfill\sf\large\leftmark}
 \vskip 1pt \hrule}}
 \def\@oddhead{\vbox{\hsize\textwidth 
 \hbox to \textwidth{{\sf\large\leftmark}
 \rule[-.6ex]{0mm}{2mm} \hfill\sf\mdseries{\thepage}}
 \vskip 1pt \hrule}}}

\def\ps@myfancyplain{
 \def\@evenhead{\vbox{\hsize\textwidth%
 \rule[-.6ex]{0mm}{2mm} \hfill }
 \vskip 1pt \hrule
 \vskip\headsep
 \vskip\textheight
 \vskip1pc
 \hbox to \textwidth{\sf\mdseries\thepage 
 \rule[-6ex]{0mm}{2mm} \hfill }}
 \def\@oddhead{\vbox{\hsize\textwidth 
 \vskip 1pt\hrule
 \vskip\headsep
 \vskip\textheight
 \vskip2pc
 \hbox to \textwidth{\hfill\rule[.4ex]{1pc}{2.5pt}
 \sf\mdseries\thepage}
}}}

\def\ps@myemptyfun{
 \def\@evenhead{\vbox{\hsize\textwidth
 \rule[-.6ex]{0mm}{2mm} \hfill }
 \vskip 1pt 
 \vskip\headsep
 \vskip\textheight
 \vskip1pc
 \hbox to \textwidth{\sf\mdseries\thepage 
 \rule[-0.6ex]{0mm}{2mm} 
 \hfill }}
 \def\@oddhead{\vbox{\hsize\textwidth 
 \vskip 1pt 
 \vskip\headsep
 \vskip\textheight
 \vskip2pc
}}}

\makeatother
\providecommand{\proofname}{Demostraci\'on.}
 {\par\noindent{\it Demostraci\'on. }\nopagebreak\normalsize}%

 {\par\noindent{\it #1. }\nopagebreak\normalsize}%
 {\hfill\linebreak[2]\hspace*{\fill}$\square$\\[-1pt]}
\makeatother

\def\sqbullet{\raise.2ex\hbox{\vrule width 3.5pt height 3.5pt}}

{}

{\em}

\newcounter{substep}
\def\thesubstep{\arabic{substep}}

{\em}

\newcounter{subsubstep}
\def\thesubsubstep{\arabic{subsubstep}}

{\em}

\numberwithin{figure}{section}

{}

\newtheoremstyle{mystyle}
  {}
  {}
  {\itshape}
  {}
  {\sf \bfseries}
  {}
{ }
  {\thmname{#1}\thmnumber{{\textcolor{blue}{\, \hspace{-1mm}#2.}}}\thmnote{ (#3)}}

\theoremstyle{mystyle}
\definecolor{royalblue(web)}{rgb}{0.25, 0.41, 0.88}
\hypersetup{
    colorlinks=true,
    linkcolor=royalblue(web),
    filecolor=magenta,      
    urlcolor=cyan,
    pdftitle={Overleaf Example},
    pdfpagemode=FullScreen,
    }

\titleformat{\section}
  {\normalfont\LARGE\bfseries}{\thesection.}{1em}{}
\titleformat{\subsection}
  {\normalfont\Large\bfseries}{\thesubsection.}{1em}{}
\titleformat{\subsubsection}
  {\normalfont\normalsize\itshape}{\thesubsubsection.}{1em}{}

\newtheorem{Teor}{Theorem}[section]

\newtheorem{Prop}[Teor]{Proposition}

\newtheorem{Defi}[Teor]{Definition}

\newtheorem{Lema}[Teor]{Lemma}

\newcommand{\R}{{\mathbb R}}

\newcommand{\mail}[1]{\small\href{mailto:#1}{#1}}

\newenvironment{Abstract}
{
\begin{center}
\textbf{Abstract}\\
\vspace{0.25cm}
\begin{minipage}{14.5cm}}
{\footnotesize
\end{minipage}
\end{center}}

\begin{document}


	\begin{center}
		{\huge {\bfseries A general theory of nonlocal elasticity based on nonlocal gradients and connections with Eringen's model}\par}
		\vspace{1cm}
		
\begin{tabular}{l@{\hskip 2cm}l} 
	{\Large José C. Bellido}{\small\textsuperscript{1}} & {\Large Guillermo García-Sáez}{\small\textsuperscript{1}} \\
	\mail{josecarlos.bellido@uclm.es} & \mail{guillermo.garciasaez@uclm.es}
\end{tabular}
\vspace{5mm}

\textsc{\textsuperscript{1}ETSII, Departamento de Matem\'aticas,\\ Universidad de Castilla-La Mancha.} \\
		Campus Universitario s/n, 13071 Ciudad Real, Spain. \\ \vspace{5mm}
\end{center}
\begin{Abstract}
We develop a general theory of nonlocal linear elasticity based on nonlocal gradients with general radial kernels. Starting from a nonlocal hyperelastic energy functional, we perform a formal linearization around the identity deformation to obtain a system of nonlocal linear elasticity equations. We establish the existence and uniqueness of weak solutions for both Dirichlet and Neumann boundary conditions, proving a general Korn-type inequality for nonlocal gradients. We show that this framework encompasses Eringen's nonlocal elasticity model as a particular case, establishing an explicit connection between the two formulations. Finally, we prove localization results demonstrating that solutions to the nonlocal problems converge to their classical local counterparts in two different regimes: as the interaction horizon vanishes and, in the fractional case, as the fractional parameter approaches one. These results provide a comprehensive and unified mathematical foundation for nonlocal elasticity theories.

\end{Abstract}

\noindent {\bf Keywords:} Nonlocal linear elasticity, nonlocal gradients, Eringen model, Korn inequality, Neumann boundary conditions for nonlocal problems, localization of nonlocal problems.


\tableofcontents


\section{Introduction}

Nonlocal models in continuum mechanics have gained increasing attention over the past two decades as powerful alternatives to classical local theories for describing materials with long-range interactions, size effects, or anomalous mechanical responses. In contrast with the classical Cauchy--Born framework---where stresses depend solely on the local deformation gradient---nonlocal theories incorporate interactions across finite or infinite horizons, thereby capturing behaviors that cannot be explained by purely local constitutive laws. Among the many approaches to nonlocality, models based on nonlocal gradients have emerged as particularly appealing due to their structural similarity with classical differential operators and their compatibility with variational formulations. These models provide a natural bridge between fractional calculus, peridynamics, and classical elasticity, and they allow for a unified treatment of nonlocal hyperelasticity and its linearized counterparts.

A central motivation for the development of nonlocal gradient theories comes from the fractional Riesz gradient introduced in \cite{ShiehSpector2015,ShiehSpector2018}, and its truncated version later proposed in connection with peridynamic-type interactions \cite{BellidoCuetoMoraCorral2023}. These operators have been shown to possess a rich functional-analytic structure, including fundamental theorems of calculus, Poincaré inequalities, and compact embeddings, which enable the formulation of well-posed nonlocal variational problems. A recent work has extended this framework to general radial kernels \cite{BellidoMoraHidde2024}, providing a broad class of nonlocal gradients that retain the essential analytical properties required for applications in continuum mechanics. In particular, the contributions \cite{BellidoCuetoMora2023,BellidoMoraHidde2024,Cueto2023,CuetoHidde2025,Kreisbeck2024} have established a comprehensive theory of nonlocal Sobolev spaces, nonlocal calculus, and $\Gamma$-convergence results for energies involving general nonlocal gradients. These developments form the mathematical foundation upon which the present work is built.

Nonlocal elasticity, both in its nonlinear and linearized forms, is a natural framework in which to apply these tools. Hyperelastic models based on nonlocal gradients allow the stored energy to depend on nonlocal deformation measures, thereby incorporating long-range interactions directly into the constitutive law. The resulting Euler--Lagrange equations involve nonlocal divergence operators and give rise to boundary-value problems that differ significantly from their classical counterparts. In the nonlinear setting, existence of minimizers typically relies on polyconvexity assumptions adapted to the nonlocal context, as established in earlier works on fractional and truncated Riesz gradients \cite{Bellido_NH2020,Bellido_NH2024}. The linearization of these models around the identity deformation leads to nonlocal analogues of the classical Navier--Lamé system, where the elasticity tensor acts on the symmetric part of a nonlocal gradient. This linearization procedure, while formally similar to the classical one, requires careful handling of the nonlocal operators and their structural properties (see \cite{BellidoCuetoMora2023}).

The linearized nonlocal elasticity model obtained in this way exhibits strong connections with Eringen's nonlocal elasticity theory, one of the most widely used nonlocal models in engineering applications. In Eringen's formulation, nonlocality enters through a stress convolution kernel, whereas in the gradient-based approach it arises through the deformation measure itself. Surprisingly, as shown in \cite{BellidoCuetoMora2023}, for the fractional and truncated Riesz gradients, certain choices of the interaction kernel in Eringen's model make the two formulations equivalent. This equivalence not only clarifies the mathematical structure underlying Eringen's theory but also provides a rigorous framework for its analysis, including existence, uniqueness, and localization results. In this paper we extend this connection to the general-kernel setting, showing that the Eringen model can be recovered from the linearized nonlocal elasticity equations associated with any admissible nonlocal gradient.

A second major theme of this work concerns localization phenomena. A desirable feature of any nonlocal model is that it should recover the classical local theory in appropriate limits. For nonlocal gradients, two natural localization regimes arise: the vanishing-horizon limit, in which the interaction radius tends to zero, and the fractional-to-local limit, in which the order of a fractional gradient tends to one. Recent results have established $\Gamma$-convergence of nonlocal energies to their local counterparts in both regimes, as well as convergence of solutions to the associated boundary-value problems \cite{BellidoCuetoMoraCorral2021,Cueto2023,CuetoHidde2025,Kreisbeck2024}. In the present work we apply these results to the linearized nonlocal elasticity model, thereby showing that the nonlocal Navier--Lamé system converges to the classical one as the horizon shrinks or, in the fractional case, as the order of fractionality approaches one. We also address the Neumann problem, where the presence of nonlocal normal derivatives introduces additional analytical challenges. While a general Neumann theory for arbitrary kernels is still under development, in the particular case of the truncated Riesz a theory was develop in \cite{Kreisbeck2024}, including localization as the fractionality parameter goes to one. Here we also address the case of the horizon going to zero.

The aim of this paper is therefore threefold. First, we perform the formal linearization of nonlocal hyperelasticity  models and establish well-posedness of the resulting nonlocal linear elasticity equations, including a general Korn-type inequality for nonlocal gradients. Both Dirichlet and Neumann boundary conditions are considered. Second, we explore the connections with Eringen's model and provide a unified perspective that encompasses fractional, truncated fractional, and general-kernel nonlocal elasticity. Finally, we present localization results for both Dirichlet and Neumann problems, relying on recent $\Gamma$-convergence theorems for nonlocal energies.

The paper is organized as follows. Section~2 reviews the necessary preliminaries on nonlocal gradients, nonlocal Sobolev spaces, and the associated calculus. Section~3 introduces the general nonlinear hyperelastic model and its linearization. Section~4 establishes existence and uniqueness of solutions for the linearized equations under both Dirichlet and Neumann boundary conditions. Section~5 develops the connection with Eringen's model and proves the equivalence of the two formulations for suitable kernels. Section~6 presents the localization results, including homogenization with respect to the horizon parameter and the fractional-to-local limit. Together, these results provide a comprehensive and unified theory of nonlocal elasticity based on general nonlocal gradients.

\section{Preliminaries}

\subsection{Notation}
We fix $n\in \mathbb{N}$ the dimension of our ambient space $\R^n$ and we will denote by $\Omega\subset \R^n$ and open bounded subset representing the body. The notation for Sobolev $W^{1,p}$ and Lebesgue $L^p$ spaces is the standard one, as is that of smooth
functions of compact support $C_c^\infty$. We will
indicate the domain of the functions, as in $L^p
(\Omega$); the target is indicated only if it is not $\R$.

We write $|x|=\left(\sum_{j=1}^n x_j^2\right)^{1/2}$ for the Euclidean norm of a vector $x=(x_1,\ldots,x_n)\in \R^n$, and $|A|$ for the Frobenius norm of a matrix $A\in \R^{m\times n}$. The identity matrix will be denoted as $I$ and the identity map as $\operatorname{id}.$ The subset of the symmetric matrices is denoted by $\R^{n\times n}_{sym}$, and the symmetric part of a matrix $M\in \R^{n\times n}$, will be denoted as $$M_{sym}=\frac{M+M^T}{2},$$ where $M^T$ is the transpose matrix of $M$. The inner product in $\R^n$ will be denoted by $\cdot$, while the respective one in $\R^n\times \R^n$ will be denoted by $:$, i.e., given $A=(a_{ij})_{1\leq i,j\leq n}$ and $B=(b_{ij})_{1\leq i,j\leq n}$, $$A:B=\sum_{i,j=1}^n a_{ij}b_{ij}.$$The tensor product of two vectors in $\R^n$ will be $\otimes$, and $*$ represents the usual operation of convolution.

The ball centered at $x\in \R^n$ and with radius $r>0$ is denoted by $B(x,r)=\{y\in \R^n: |x-y|<r\}$. The complementary of a set $E\subset \R^n$ is denoted by $E^c=\R^n\setminus{E}$, its closure by $\overline{E}$ and its boundary by $\partial E$. The distance between a point $x\in \R^n$ and an open set $U\subset\R^n$ is denoted by $d(x,E)$. For $\delta>0$ we denote $$\Omega_\delta:=\Omega+B(0,\delta)=\{x\in \R^n: d(x,\Omega)<\delta\},$$ and $$\Omega_{-\delta}=\{x\in \Omega: d(x,\Omega^c)>\delta\},$$ whenever it makes sense. In the case that $\Omega_{-\delta}$ is well defined, we introduce the \textit{nonlocal boundary}, denoted by $\Gamma_{\pm \delta}$, as $$\Gamma_{\pm \delta}:=\Omega_{\delta}\setminus{\Omega_{-\delta}}=\Gamma_\delta\cup \Gamma_{-\delta}\cup\partial\Omega,$$ where $$\Gamma_\delta:=\Omega_\delta\setminus\overline{\Omega},\,\Gamma_{-\delta}:=\Omega\setminus{\overline{\Omega_{-\delta}}},$$ are the inner and outer collars of $\Omega$, respectively.

Our convention for the Fourier transform of functions $f\in L^1(\R^n)$ is $$\widehat{f}(\xi)=\int_{\R^n}f(x)e^{-2\pi x\cdot \xi}\,dx,\,\xi\in \R^n.$$ This definition is extended by continuity and duality to other function and distribution spaces as usually in function spaces theory. We will sometimes use the alternative notation $\mathfrak{F}(f)$ for $\widehat{f}$. More details of this operator could be found in the classical texts \cite{Duoandikoetxea2001,Grafakos2008}.

Regarding radial functions we have the following definitions:
\begin{itemize}
    \item A function $f:\R^n\to \R$ is \textit{radial} if there exists $\overline{f}:[0,\infty)\to \R$ such that $\overline{f}(|x|)=f(x)$ for every $x\in \R^n$. The function $\overline{f}$ is called the \textit{radial representation} of $f$.
    \item A radial function $f:\R^n\to \R$ is \textit{radially decreasing} if its radial representation is a
decreasing function.
\item A function $g:\R^n\to \R^n$ is \textit{vector radial} if there exists a radial function $\overline{\phi}:[0,\infty)\to \R$  such that
$\phi(x)=\overline{\phi}(|x|)x$ for every $x\in \R^n$.
\end{itemize}

For real valued functions, we use the monotonicity properties of being increasing and decreasing in
the non-strict sense. A function $f :\R \to\R$ is called \textit{almost decreasing} if there exists a positive constant $C$ such
that $f(t)\geq Cf(s)$ for every $t\leq s$, and an analogous definition for \textit{almost increasing}.

A function $f:\Omega\times \R^m\times \R^{m\times n}$ is called a \textit{Carath\'eodory integrand} if $f(\cdot, z,A)$ is measurable for all $(z,A)\in \R^m\times \R^{m\times n}$ and $f(x,\cdot,\cdot)$ is continuous a.e. $x\in \Omega$.

\subsection{Nonlocal gradients}

In this section, we introduce the key elements of our functional setting, namely the nonlocal gradients for general kernels and their associated function spaces. Nonlocal gradients for general kernels were introduced in \cite{MengeshaSpector}, and the calculus for those nonlocal objects has been addressed in \cite{DEliaGulianMengeshaScott2022,DEliaGulianOlsonKarniadakis2021,DuGunzburgerLehoucqZhou2013}. Of particular interest to us are the functional spaces associated with these objects, and in particular, structural properties such as Poincar\'e inequalities and compact embeddings. This has been developed for general nonlocal gradients with radial kernels in \cite{BellidoMoraHidde2024}, generalizing the case of the {\it truncated} Riesz fractional gradient analyzed in \cite{BellidoCuetoMoraCorral2023} (see also \cite{CuetoHidde2025}). The case of non-radial kernels has also been considered in \cite{HanTian2023}.

\subsection{Definitions and nonlocal calculus}

Throughout the entire manuscript, we will assume that the kernel in the definition of nonlocal gradients $\rho$ is a function that verifies:
\begin{equation}\label{eqn:H0} \tag{H0}\begin{cases} \text{$\rho:\R^n\setminus\{0\} \to [0,\infty)$ is radial}, \\
\text{$\rho \in L^1_{\text{loc}}(\R^n)$ with $\displaystyle \int_{\R^n} \min\{1,\abs{x}^{-1}\}\rho(x)\,dx < \infty$,} \\
\text{$\inf_{\overline{B(0,\varepsilon)}}\rho >0$ for some $\varepsilon>0$.} \end{cases}\end{equation}
Hypotheses, as (\ref{eqn:H0}), are somehow the minimal requirement for the kernel $\rho$ in the radial case and have frequently appeared in the related literature. Observe that, since by \eqref{eqn:H0} the kernel $\rho$ is radial, there exists a function $\overline{\rho}:(0,+\infty)\to [0,+\infty)$ such that $\rho(x)=\overline{\rho}(|x|)$. With a slight abuse of notation, and without risk of ambiguity, we shall henceforth denote $\overline{\rho}$ simply by $\rho$.

\begin{Defi}[Nonlocal gradient]\em
    Let $u\in C_c^\infty(\R^n)$, we define the \textit{nonlocal gradient} with kernel $\rho$ of $u$ as $$D_{\rho}u(x):=\int_{\R^n}\frac{u(x)-u(y)}{|x-y|}\frac{x-y}{|x-y|}\rho(x-y)\,dy,\,x\in \R^n,$$ and for $u\in C_c^\infty(\R^n;\R^n)$ as $$D_\rho u(x):=\int_{\R^n}\frac{u(x)-u(y)}{|x-y|}\otimes \frac{x-y}{|x-y|}\rho(x-y)\,dy,\,x\in \R^n.$$
\end{Defi}

It is straightforward to see that for functions $u\in C_c^\infty(\R^n)$, the integral involving the definition of $D_\rho u$ is absolutely convergent for each $x\in \R^n$. Moreover, $D_\rho u\in L^1(\R^n;\R^n)\cap L^\infty(\R^n;\R^n)$ (\cite{BellidoMoraHidde2024}), and hence $D_\rho u\in L^p(\R^n;\R^n)$ for each $1<p<\infty$ by Riesz-Torin interpolation inequality.\\

In \cite{BellidoMoraHidde2024, CuetoHidde2025}, examples of admissible kernels are given, the most paradigmatic one 
being $$\rho^s(x):=c_{n,s}\frac{1}{|x|^{n+s-1}},\,s\in (0,1),$$ where $c_{n,s}$ is a suitable normalization constant. The nonlocal gradient associated with this kernel is the well-known Riesz fractional gradient $D^s:=D_{\rho^s}$, which is defined for functions $u\in C_c^\infty(\R^n)$ as $$D^su(x)=c_{n,s}\int_{\R^n}\frac{u(x)-u(y)}{|x-y|^{n+s}}\frac{x-y}{|x-y|}\,dy.$$
 The Riesz fractional gradient was introduced by Shieh and Spector in \cite{ShiehSpector2015,ShiehSpector2018}, and since then has been widely studied.\\
 
Another important example concerns a truncated version of the Riesz fractional gradient introduced in \cite{BellidoCuetoMoraCorral2023}, which is suitable for bounded domains. This formulation is motivated by concepts from peridynamics, as its integration domain is restricted to a ball of radius $\delta>0$, corresponding to the interaction horizon among particles in the peridynamic framework. Given $\delta>0$, and a cut-off function $w_\delta:\mathbb{R}^n\to [0,+\infty)$ verifying
\begin{itemize}
    \item $w_\delta$ is radial;
    \item $w_\delta\in C_c^\infty\left(B(0,\delta)\right)$;
    \item There exist constants $a_0>0$ and $b_0\in (0,1)$ such that $0\leq w_\delta\leq a_0$, and $w_\delta|_{B(0,b_0\delta)}=a_0$.
    \item $w_\delta$ is radially decreasing.
    \item The function $w_\delta$ is normalized in the following sense: $$\int_{B(0,\delta)}\frac{w_\delta(x)}{|x|^{n+s-1}}=\frac{n}{c_{n,s}}.$$ 
\end{itemize}
Now, taking the kernel 
$$\rho_\delta^s(x):=c_{n,s}\frac{w_\delta(x)}{|x|^{n+s-1}},$$ 
yields that $D_{\rho_\delta^s}=D_\delta^s$, where $D_\delta^s$ is the nonlocal gradient introduced in \cite{BellidoCuetoMoraCorral2023}.\\

The theory established for $D^s$ and $D^s_\delta$ has been further generalized to encompass nonlocal gradients $D_\rho$, where the kernel $\rho$ satisfies hypothesis \eqref{eqn:H0}, as shown in \cite{BellidoMoraHidde2024}. The present work is developed within this broader framework, and in this section we provide an overview of the fundamental aspects of the calculus involved, together with the underlying functional spaces.

 The first important result is the following representation formula for $D_\rho$ in terms of a convolution with the locally integrable function 
 $$Q_\rho(x):=\int_{|x|}^\infty \frac{\overline{\rho}(t)}{t}\,dt,\,x\in \R^n\setminus\{0\}.$$
 By \cite[Proposition 2.6]{BellidoMoraHidde2024}, we have that $$D_\rho u=Q_\rho*Du=D\left(Q_\rho*u\right),$$ 
 for $u\in C_c^\infty(\R^n)$, and moreover, if $\rho\in L^1(\R^n)$, then $Q_\rho\in L^1(\R^n)$ (in fact, $\rho$ having compact support is enough for $Q_\rho$ to lie in $L^1(\R^n)$, \cite[Lemma 2.5 (iii)]{BellidoMoraHidde2024}). Also, taking Fourier transform, 

 $$\widehat{D_\rho u}(\xi)=2\pi i\xi \widehat{Q_\rho}(\xi)\widehat{u}(\xi),$$ with $$\widehat{Q_\rho}(\xi)=\frac{1}{2\pi|\xi|}\int_{\R^n}\frac{\rho(x)x_1}{|x|^2}\sin(2\pi|\xi|x_1)\,dx,\,\xi\in \R^n\setminus{\{0\}}.$$

Note that the function $Q_\rho$ plays the role of the Riesz potential $$I_{1-s}(x):=\frac{1}{\gamma_{1-s}}\frac{1}{|x|^{n-s+1}},\,0<1-s<n,$$ where $$\gamma_{1-s}=\frac{n+s-1}{c_{n,s}},$$ for the Riesz fractional gradient $D^s$. In fact, for the choice $\rho(x)=c_{n,s}\frac{1}{|x|^{n+s-1}}$, we have that \begin{align*}
    Q_\rho(x)=c_{n,s}\int_{|x|}^\infty \frac{t^{-n-s+1}}{t},\,dt=c_{n,s}\int_{|x|}^\infty t^{-(n+s)}\,dt=\frac{c_{n,s}}{n+s-1}\frac{1}{|x|^{n+s-1}}=I_{1-s},
\end{align*}
and hence $D^su=D(I_{1-s}*u)$, as it was proven in \cite[Theorem 1.2]{ShiehSpector2015}. A similar reasoning shows that for the choice $\rho(x)=c_{n,s}w_\delta(x)\frac{1}{|x|^{n-1+s}}$ yields that $Q_\rho=Q_\delta^s(x)$, where $Q_\delta^s$ is the function of \cite[Proposition 4.3]{BellidoCuetoMora2023}.\\

We can also define the nonlocal divergence for general kernels $\rho$.
\begin{Defi}\em
  For $v\in C_c^\infty(\R^n;\R^n)$, we define the {\it nonlocal divergence} with kernel $\rho$ as 
  $$\operatorname{div}_\rho v(x):=\int_{\R^n}\frac{v(x)-v(y)}{|x-y|}\cdot \frac{x-y}{|x-y|}\rho(x-y)\,dy,\,x\in \R^n,\,$$ and for $v\in C_c^\infty(\R^n,\R^{n\times n})$ as  $$\operatorname{div}_\rho v(x):=\int_{\R^n}\frac{v(x)-v(y)}{|x-y|} \frac{x-y}{|x-y|}\rho(x-y)\,dy,\,x\in \R^n,.$$
\end{Defi}
In this case we also have a duality relationship between $D_\rho$ and $\operatorname{div}_\rho$ as in the fractional and finite horizon cases by \cite[Proposition 3.2]{BellidoMoraHidde2024}.
\begin{Lema}[Integration by parts]\label{IBP}
    Let $u\in C_c^\infty(\R^n)$ and $v\in C_c^\infty(\R^n;\R^n)$. Then, $$\int_{\R^n}D_\rho u\cdot v\,dx=-\int_{\R^n}u\operatorname{div}_\rho v\,dx.$$
\end{Lema}Again, by the suitable choice for $\rho$, we can recover the fractional divergence $\operatorname{div}^s$ and the nonlocal counterpart $\operatorname{div}_\delta^s$, and hence the integration by parts formulas of the fractional case \cite[Lemma 2.2]{BellidoCuetoMoraCorral2021} and the finite horizon one \cite[Theorem 3.2]{BellidoCuetoMoraCorral2023}.\\

\subsection{Nonlocal Sobolev spaces}

In analogy to the fractional Sobolev spaces $H^{s,p}(\R^n)$, for $1<p<\infty$, the $\rho$-nonlocal Sobolev spaces $H^{\rho,p}(\R^n)$ were introduced in \cite{BellidoMoraHidde2024} as $$H^{\rho,p}(\R^n):=\{u\in L^p(\R^n): D_\rho u\in L^p(\R^n;\R^n)\},$$ endowed with the norm $$\norm{u}_{\rho,p}:=\norm{u}_p+\norm{D_\rho u}_p.$$ As a consequence of \cite[Theorem 3.9 (i)]{BellidoMoraHidde2024} and the completeness of those spaces, we can alternatively define them as $$H^{\rho,p}(\R^n):=\overline{C_c^\infty(\R^n)}^{\norm{\cdot}_{\rho,p}}.$$
For the choice $\rho=\rho^s, s\in (0,1)$, we have that by \cite[Theorem 1.7]{ShiehSpector2015}, $$H^{\rho,p}(\R^n)=H^{s,p}(\R^n),$$ where $H^{s,p}(\R^n)$ is the classical Bessel potential space (see \cite{BellidoGarcia2025, BellidoCuetoGarcia2025}). 
For an open set $\Omega\subset\R^n$, we define the closed subspace $$H^{\rho,p}_0(\Omega):=\left\{u\in H^{\rho,p}(\R^n): u(x)=0\text{ a.e. }x\in\Omega^c\right\}.$$ 
If $\Omega$ has a Lipschitz domain, by \cite[Theorem 3.9 (iii)]{BellidoMoraHidde2024} we have that $$H^{\rho,p}_0(\Omega)=\overline{C_c^\infty(\Omega)}^{\norm{\cdot}_{\rho,p}},
$$ where the elements $C_c^\infty(\Omega)$ should be interpreted as its extension to $\R^n$ by zero. Since for $p=2$, the spaces $H^{\rho,p}(\R^n)$ and  $H^{\rho,p}_0(\Omega)$ are Hilbert, we denote them by $H^\rho(\R^n)$ and $H^\rho_0(\Omega)$, in analogy to the classical Sobolev theory.

In order to ensure the fundamental structural properties of the spaces $H^{\rho,p}$, such as Poincar\'e-type inequalities and compact embeddings, it is necessary to impose additional assumptions on the kernels $\rho$. The first one is that we will restrict ourselves to kernels with compact support, which can be assumed without loss of generality since the main properties of the function spaces $H^{\rho,p}(\R^n)$ and $H_0^{\rho,p}(\Omega)$ are encoded on the behavior of $\rho$ near
zero by \cite[Proposition 3.10]{BellidoMoraHidde2024}.

Let $\varepsilon$ be as in \eqref{eqn:H0}. We will employ the following conditions:
\begin{enumerate}
    \item[(H1)]\label{eqn:H1} There exists $\nu>0$ such that the functions $f_\rho:(0,\infty)\to \R$, defined by $f_\rho(r)= r^{n-2}{\rho}(r)$, and $g(r)= r^\nu f_\rho(r)$ are decreasing on $(0,\varepsilon)$;
    \item[(H2)]\label{eqn:H2} The function $f_\rho$ is smooth away from the origin, and for every positive integer $k$ there exists a constant $C(k)>0$ such that

\[
    \left|\frac{d^kf_\rho(r)}{dr^{k}}\right|\leq C(k)\,r^{-k}f_\rho(r), \quad r\in (0,\varepsilon);
    \]

    \item[(H3)]\label{eqn:H3} There exists $s\in(0,1)$ such that the function $r\mapsto r^{n+s-1}{\rho}(r)$ is almost decreasing on $(0,\varepsilon)$;
    \item[(H4)]\label{eqn:H4} There exists $t\in(0,1)$ such that the function $r\mapsto r^{n+t-1}{\rho}(r)$ is almost increasing on $(0,\varepsilon)$.
\end{enumerate}

These hypotheses are natural, and we refer to Sections~4 and~5 of \cite{BellidoMoraHidde2024} for a detailed discussion. As expected, they are satisfied by the kernels associated with $D^s$ and $D_\delta^s$. Assuming $(H0)$–$(H4)$ for $\rho$, we obtain the following result \cite[Theorem~5.2]{BellidoMoraHidde2024}:

\begin{Teor}[General nonlocal fundamental theorem of calculus]
    Let $\rho$ have compact support and satisfy $(H0)$-$(H4)$. Then there exists a radial vector function

\[
    V_\rho \in C^\infty\!\left(\R^n\setminus\{0\};\R^n\right)\cap L_{\text{loc}}^1(\R^n;\R^n)
    \]

    such that for all $u\in C_c^\infty(\R^n)$,

\[
    u(x)=V_\rho*D_\rho u(x), \quad x\in \R^n.
    \]

    Moreover, there exists a constant $C>0$ such that for all $x\in B_\varepsilon(0)\setminus\{0\}$, with $\varepsilon$ as in \eqref{eqn:H0},

\[
    |V_\rho(x)|\leq \frac{C}{|x|^{2n-1}\rho(x)}, \qquad |DV_\rho(x)|\leq \frac{C}{|x|^{2n}\rho(x)}.
    \]

\end{Teor}

Note that if we consider the kernel of the nonlocal gradient $D_\delta^s$, we recover the nonlocal fundamental theorem of calculus of \cite[Theorem 4.5]{BellidoCuetoMoraCorral2023}.\\

Under the same hypothesis, Poincar\'e and compactness results were proven for general kernels with suitable behaviour around the origin. In particular, Poincar\'e inequality will be a key result for our aim in this paper.
\begin{Teor}\label{Poincare}
    Let $p\in (1,\infty)$, $\Omega\subset \R^n$ an open bounded set and $\rho$ satisfying (H0)-(H2) and having compact support. Then, the following two statements hold:
    \begin{enumerate}
        \item (Poincar\'e inequality) If $$\liminf_{t\to 0^+}t^{n-1}\overline{\rho}(t)>0,$$ then there is a positive constant $C$ such that $$\norm{u}_{L^p(\Omega)}\leq C\norm{D_\rho u}_{L^p(\R^n;\R^n)}\,$$
        for all $u\in H_0^{\rho,p}(\Omega)$
        \item (Compact embedding) If $$\lim_{t\to 0^+}t^{n-1}\overline{\rho}(t)=\infty,$$ then $H^{\rho,p}_0(\Omega)$ is compactly embedded into $L^p(\R^n)$.
    \end{enumerate}
\end{Teor}
Many different embeddings in the setting of Orlicz spaces were proven in \cite[Theorem 6.2]{BellidoMoraHidde2024}. In particular, we are interested in the ones on the Lebesgue scale.
\begin{Teor}\label{Embedding}
    Let $\rho$ satisfy (H1)-(H4) with compact support and pick $s,t$ as in those hypothesis. Let $1<p<\infty$ such that $tp<n$. Then, $$H^{\rho,p}_0(\Omega)\xhookrightarrow{}L^q(\Omega),\,1\leq q\leq p_s^*,$$ where $p_s^*=\frac{np}{n-sp}$.
\end{Teor}
Moreover, when $sp>n$, \cite[Theorem 6.5]{BellidoMoraHidde2024}, yields that the embedding extends to the case $q=\infty$.\\

The following lemma, obtained by extending the representation formula for $D_\rho$ in terms of convolution with $Q_\rho$ via density arguments (see \cite[Lemma~2.9]{CuetoHidde2025}), provides a means to translate general nonlocal gradients into their classical counterparts. It will serve as a key tool in establishing Korn's inequality for $D_\rho$. Since a proof is not included in \cite[Lemma~2.9]{CuetoHidde2025}, we present here a complete argument, in the spirit of \cite[Theorem~3]{Cueto2023}, for the reader's convenience.

\begin{Lema}[From nonlocal to local]\label{NLTOL}
    Let $\rho$ satisfy (H1)-(H4) with compact support. The linear map $\mathcal{Q}_\rho$ defined as $\mathcal{Q}_\rho u:=Q_\rho *u$, $u\in C_c^\infty(\R^n)$, extends to a linear bounded map from $H^{\rho,p}(\R^n)\to W^{1,p}(\R^n)$. Furthermore, it holds for all $u\in H^{\rho,p}(\R^n)$ that $$D_\rho u=D\left(\mathcal{Q}_\rho u\right).$$
\end{Lema}
\noindent{}\textbf{Proof:} Let $u\in H^{\rho,p}(\R^n)$. Since $\rho$ has compact support, $Q_\rho\in L^1(\R^n)$, and by Young's inequality for convolutions we have that $v:=\mathcal{Q}_\rho u\in L^p(\R^n)$, hence $\mathcal{Q}_\rho$ is a bounded linear operator. Now, for every $w\in C_c^\infty(\R^n;\R^n)$, \begin{align*}\int_{\R^n}v\operatorname{div}w\,dx&=\int_{\R^n}(Q_\rho*u)\operatorname{div}w\,dx=\int_{\R^n}u(Q_\rho*\operatorname{div}w)\,dx=\int_{\R^n}u\operatorname{div}_\rho w\,dx\\
&=-\int_{\R^n}D_\rho u\cdot w\,dx.
\end{align*} where we have used the Fubini's theorem for the first identity and Lemma \ref{IBP} for the last one. By the distributional definition of $W^{1,p}(\R^n)$, this proves that $v\in W^{1,p}(\R^n)$ with $Dv=D_\rho u$, and hence $\mathcal{Q}_\rho u\in W^{1,p}(\R^n)$ for every $u\in H^{\rho,p}(\R^n)$.\qed

In \cite[Lemma~2.12]{CuetoHidde2025}, an inverse result to the previous one is established under the same assumptions on the kernel $\rho$. Unlike the earlier case, this lemma includes the proof and shows the existence of an inverse translation operator from $W^{1,p}$ to $H^{\rho,p}$.

\begin{Lema}[From local to nonlocal]\label{LNL} Let $\rho$ satisfy $(H1)-(H4)$ with compact support. The linear map $\mathcal{P}_\rho$ defined as $$v\mapsto\mathcal{P}_\rho  v:=\mathfrak{F}^{-1}\left(\widehat{v}/\widehat{Q}_\rho\right),\,v\in \mathcal{S}(\R^n),$$ extends to a bounded operator from $W^{1,p}(\R^n)$ to $H^{\rho,p}(\R^n)$. Moreover, $(\mathcal{P}_\rho)^{-1}=\mathcal{Q}_\rho$, and for every $v\in W^{1,p}(\R^n)$, $$Dv=D_\rho(\mathcal{P}_\rho v).$$
\end{Lema}
Note that Lemmas \ref{NLTOL} and \ref{LNL} are also valid for the spaces $W^{1,p}_0(\Omega)$ and $H^{\rho,p}_0(\Omega)$, whenever $\Omega$ has Lipschitz boundary, via extension by zero outside $\Omega$ to the whole $\R^n$.\\
As it was pointed in \cite[Remark 4(e)]{Cueto2023} for the nonlocal gradient $D_\delta^s$, Lemmas \ref{NLTOL} and \ref{LNL} gives us an alternative way to prove Poincar\'e inequality under slightly different, but in the end equivalent, conditions as the ones of Theorem \ref{Poincare}.
\begin{Teor}[Poincar\'e inequality]\label{PoincareII}
    Let $\rho$ a kernel satisfying $(H0)-(H4)$ with compact support. Then, there exists a constant $c>0$, such that $$\norm{u}_{L^p(\Omega)}\leq c\norm{D_\rho u}_{L^p(\R^n;\R^n)},\,u\in H_0^{\rho,p}(\Omega).$$
\end{Teor}
\noindent{}\textbf{Proof:} Let $u\in H_0^{\rho,p}(\Omega)$. Since by Lemmas \ref{NLTOL} and \ref{LNL}, $\mathcal{Q}_\rho:H^{\rho,p}_0(\Omega)\to W^{1,p}_0(\Omega)$ and $\mathcal{P}_\rho:W^{1,p}_0(\Omega)\to H^{\rho,p}_0(\Omega)$, continuously, with $(\mathcal{P}_\rho)^{-1}=\mathcal{Q}_\rho$, the classical Poincar\'e inequality implies that \begin{align*}
    \norm{u}_{L^p(\Omega)}=\norm{\mathcal{P}_\rho \mathcal{Q}_\rho u}_{L^p(\R^n)}\leq c\norm{\mathcal{Q}_\rho u}_{W^{1,p}_0(\R^n)}\leq c\norm{D\mathcal{Q}_\rho u}_{L^p(\R^n;\R^n)}=c\norm{D_\rho u}_{L^p(\R^n;\R^n)},
\end{align*}
as we wanted to prove. \qed

\subsection{Poincar\'e-Wirtinger inequality and Neumann boundary conditions}

Given a bounded open set $\Omega \subset \R^n$, we established in the previous section a Poincar\'e inequality within the closed subspace $H_0^{\rho,p}(\Omega)$, namely the subspace of $H^{\rho,p}(\R^n)$ consisting of functions vanishing on the complement of $\Omega$. In analogy with the classical Sobolev spaces, it is of significant interest to derive related inequalities when this Dirichlet condition on the complement is relaxed. In this section, we address this problem following \cite{Kreisbeck2024}, which provides Poincar\'e--Wintiger type inequalities for the truncated Riesz fractional gradient $D^s_\delta$.

We start with an arbitrary kernel $\rho$ satisfying $(H0)$--$(H4)$, such that $\operatorname{supp}\rho=\overline{B(0,1)}$. 

Given a horizon $\delta>0$ such that $\Omega_{-\delta}$ is well defined, the rescaled kernel $\rho_\delta$ is defined as
\[
\rho_\delta(x):=\delta^{-n}\rho(x/\delta), \quad x\in \R^n\setminus\{0\}.
\]
The choice of the scaling factor implies that
\[
\widehat{Q_{\rho_\delta}}(\xi)=\widehat{Q_\rho}(\delta\xi).
\]

The nonlocal Sobolev spaces associated with $D_{\rho_\delta}$ are defined analogously, but we denote them as
\[
H^{\rho,p,\delta}(\R^n)=H^{\rho_\delta,p}(\R^n), \qquad H^{\rho,p,\delta}_0(\Omega)=H^{\rho_\delta,p}_0(\Omega),
\]
to emphasize the dependence on the parameter $\delta$. As in the case of the nonlocal gradient $D_\delta^s$ (see \cite{BellidoCuetoMoraCorral2023}), the definition of $D_{\rho_\delta}u$ on $\Omega$ depends only on the values of $u$ in $\Omega_\delta$. Typically,  homogeneous Dirichlet conditions are $u=0$ on $\Gamma_{\pm\delta}$.

In order to obtain Poincar\'e--Wirtinger type inequalities, it is critical to understand the set of functions $u\in H^{\rho,p,\delta}(\Omega)$ such that $D_{\rho_\delta}u=0$. This is a challenging problem that has been successfully addressed in \cite{Kreisbeck2024} in the particular case $D_{\rho}=D_\delta^s$. In this section we focus on this specific case.

In \cite[Proposition 3.1]{Kreisbeck2024} it was shown that the space
\[
N^{s,p,\delta}(\Omega):=\{u\in H^{s,p,\delta}(\Omega):D^s_\delta u=0\},
\]
contains functions that are not constant a.e. in $\Omega$. Moreover, in \cite[Proposition 3.3]{Kreisbeck2024} it was proven that $N^{s,p,\delta}(\Omega)$ is an infinite-dimensional closed subspace of $H^{s,p,\delta}(\Omega)$. In particular, this implies that no Poincar\'e inequality is available for the space $H^{s,p,\delta}(\Omega)$, hence we must seek a suitable subspace for our model. To this end, the nonlocal metric projection $\pi_\delta^s:L^p(\Omega_\delta)\to N^{s,p,\delta}(\Omega)$ is introduced as the function minimizing the distance to functions with vanishing nonlocal gradient in the $L^p$-norm, i.e.,
\[
\norm{u-\pi_\delta^s(u)}_{L^p(\Omega_\delta)}=\min_{v\in N^{s,p,\delta}(\Omega)}\norm{u-v}_{L^p(\Omega_\delta)}.
\]

From this, we define
\[
N^{s,p,\delta}(\Omega)^\perp:=\{u\in H^{s,p,\delta}(\Omega):\pi_\delta^s(u)=0\}.
\]

In the case $p=2$, the function $\pi_\delta^s$ is the linear orthogonal projection onto $N^{s,2,\delta}(\Omega)$, and hence $N^{s,2,\delta}(\Omega)^\perp$ is the orthogonal complement of $N^{s,2,\delta}(\Omega)$. For $p\neq 2$, $\pi_\delta^s$ need not be linear, but it is continuous for every $p$.

The relevant mathematical fact about this projection is that it allows us to establish a nonlocal Poincar\'e--Wirtinger inequality, as shown in \cite[Lemma 4.7]{Kreisbeck2024}.
\begin{Teor}\label{PW}
    Let $1<p<\infty$. Then, there exists a constant $C$ independent of $s$ such that

\[
    \norm{u-\pi_\delta^s(u)}_{L^p(\Omega_\delta)}\leq C\norm{D^s_\delta u}_{L^p(\Omega;\R^n)},
    \]

    for every $u\in H^{s,p,\delta}(\Omega)$.
\end{Teor}

We finally recall the following nonlocal operator introduced in \cite{BellidoCueto2024}:
\begin{Defi}\em
    Let $u:\Omega_\delta\to \R$ and $v:\Omega_\delta\to \R^n$ be two measurable functions. We define the \textit{nonlocal normal derivative} at $x\in \Gamma_{\pm\delta}$ of $u$ and $v$ as
    \begin{align*}
    \mathcal{N}_\delta^su(x)&:=c_{n,s}\operatorname{pv}_x\int_\Omega\frac{u(x)-u(y)}{|x-y|^{n+s}}\frac{x-y}{|x-y|}w_\delta(x-y)\,dy,\\
    \mathcal{N}_\delta^sv(x)&:=c_{n,s}\operatorname{pv}_x\int_\Omega\frac{v(x)-v(y)}{|x-y|^{n+s}}\cdot\frac{x-y}{|x-y|}w_\delta(x-y)\,dy,
    \end{align*}
    whenever this makes sense.
\end{Defi}
This operator arises naturally in the integration by parts formula for functions $u\in H^{s,p,\delta}(\Omega)$ for which boundary conditions are not prescribed.

\section{A general model via linearization}

In this section we perform the formal linearization of nonlinear hyperelasticity models. We follow the lines of \cite[Section 5]{BellidoCuetoMora2023}.
From now on, we will assume that the kernel $\rho$ is in the hypothesis of Theorem \ref{PoincareII}.
Moreover, we assume that the kernel is normalized in the following way 
$$\int_{\R^n}\rho (x)\,dx=n.$$ 
We require this in order to have the following technical result, which establishes that the nonlocal gradient of an affine map is the associated matrix, generalizing \cite[Proposition 4.1]{BellidoCuetoMora2023}. The proof is a straightforward generalization of \cite[Proposition 4.1]{BellidoCuetoMora2023}.
\begin{Lema}
    Let $A\in \R^{n\times n}$ and $b\in \R^n$. Let $u\in H_0^{\rho,p}(\Omega;\R^n)$ defined as $u(x):=Ax+b$, $x\in \Omega$. Then, $D_\rho u(x)=A$, $x\in \Omega$.
\end{Lema}

We now introduce the following hyperelastic nonlinear energy associated with a given deformation mapping 
$u \in H^{\rho,p}_{\operatorname{id}}(\Omega;\R^n)$:
\[
I[u]=\int_\Omega W\bigl(x,D_\rho u\bigr)\,dx - \int_\Omega f \cdot u\,dx,
\]
where the integrand $W:\Omega \times \R^{n\times n} \to \R \cup \{+\infty\}$ is a Carath\'eodory function representing the stored elastic energy density, and $f:\Omega \to \R^n$ denotes a body force acting on the domain. For simplicity, we impose Dirichlet boundary conditions on the complementary domain, namely $u|_{\Omega^c}=\operatorname{id}$.

Let us assume that the integrand $W$ satisfies the conditions in \cite[Theorem 8.1, Theorem 8.2]{BellidoCuetoMoraCorral2023}, which includes as a central requirement polyconvexity of $W(x,\cdot)$ almost every $x$, and ensure the existence of minimizers for $I$, and allows us to derive under the integral sign in order to obtain the Euler-Lagrange equations. Thus, given $v\in C_c^\infty(\Omega)$, we have  
\begin{align*}
    \frac{d}{dt}I[u+tv]\Bigg|_{t=0}=\int_\Omega W_2(x,D_\rho u)\cdot D_\rho v\,dx-\int_{\Omega}f\cdot v\,dx,
\end{align*} where $W_2$ denotes the derivative of $W(x,\cdot)$. By Lemma \ref{IBP}, $$\int_\Omega W_2(x,D_\rho u)\cdot D_\rho v\,dx=-\int_{\Omega}\operatorname{div}_\rho\left(W_2(x,D_\rho)\right)\,v\,dx,$$ and hence $$-\int_{\Omega}\operatorname{div}_\rho\left(W_2(x,D_\rho u)\right)\,v\,dx=\int_\Omega f\cdot v\,dx,\,v\in C_c^\infty(\Omega),$$ and in virtue of the fundamental lemma of the calculus of variations we obtain the Dirichlet problem $$\begin{cases}
    -\operatorname{div}_\rho W_2(x,D_\rho u)&=f,\,\text{in $\Omega$},\\
    \hfill u&=\operatorname{id},\,\text{on $\Omega^c$}
\end{cases}$$
where $W_2(x,F)=T_R(x,F)$ is the Piola-Kirchhoff stress tensor, which depends on nonlinear operators as the deformation jacobian. In order to obtain simpler linear models, in the classical case, i.e., when the Euler-Lagrange equations reads as 
$$-\operatorname{div}W_2(x,Du)=f,\,x\in \Omega,$$
a linearization process is performed (see, e.g., \cite[Ch. X]{Gurtin1981}, \cite[Ch. 52]{Gurtin2010}, \cite[Ch. 4]{Marsden1994}). This process being completely formal is interesting as we can recover the classical linear Navier-Lam\'e system. We waould like to point out that a rigorous derivation of the linearization process would require $\Gamma$-convergence arguments adapting the ones in \cite{DalMaso2002}, including nonlocal versions of the rigidity estimates in \cite{Muller2002}. 

Since the nonlocal operator $\operatorname{div}_\rho$ is linear, the linearization process concerns only to the term $T_R(x,D_\rho u)$, however, in view of how the process is performed it is enough to substitute $D$ by $D_\rho$ in the classical theory. The classical theory process consists in linearizing $T_R$ around the identity map assuming that $W$ does not depend on $x$. Also, the identity is assumed to be stress-free, i.e., $T_R(I)=0$, which is a very natural condition from the mechanical point of view. If we write $u=\operatorname{id}+v$, with $v\in H_0^{\rho,p}(\Omega;\R^n)$ the displacement, we have that 
$$D_\rho u=D_\rho (\operatorname{id}+v)=D_\rho\operatorname{id}+D_\rho v=I+D_\rho v.$$
Hence, following the steps in \cite[Ch. X]{Gurtin1981}, we have that 
$$T_R(D_\rho u)=T_R(I+D_\rho v)=T_R(I)+DT_R(I)D_\rho v+o(D_\rho v).$$ 
We may assume, as in the classical case, that $D_\rho v$ is small in the sense that $o(D_\rho v)\approx0$. The tensor that arises as the first derivative of $T_R$ evaluated in the identity $I$ is called the \textit{elasticity tensor}, and it is denoted as $C:=DT_R(I)$. Since $T_R(F)=W_2(F)=DW(F)$, we have that the components of the tensor are 
$$c_{ijkl}=\frac{\partial^2 W}{\partial F_{ij}\partial F_{kl}}(I),\quad 1\leq i,j,k,l\leq n.$$
Since $C$ is a fourth order tensor acting over matrices $M$, we will denote the action of $C$ over $M$ as $C[M]$ instead of $CM$ to differentiate it from the usual product of matrices. One of the main properties of the elasticity tensor is that $C[M]=C[M_{sym}]\in \R^{n\times n}_{sym}$, hence we can consider $C$ as a map from $\R_{sym}^{n\times n}\to \R^{n\times n}_{sym}$.

Therefore, the the first-order linear of $T_R(D_\rho u)$ is $CD_\rho v$, and as in the classical case, $C[D_\rho] v=C[D_{\rho,sym} v]$, where $D_{\rho,sym}v=(D_\rho v)_{sym}$. 


Linearization yields the strong form of the equations of general nonlocal linear elasticity: \begin{equation}\label{eqn:P}\begin{cases}
    -\operatorname{div}_\rho \left(C[D_{\rho,sym}v]\right)&=f,\,\text{in $\Omega$},\\
    \hfill v&=0,\,\text{on $\Omega^c$}.
\end{cases}\end{equation}
Multiplying both sides by a function $w\in H^{\rho}_0(\Omega;\R^n)$, and using again the duality relationships of Lemma \ref{IBP}, we get that the weak form of the equations is the following: find $v\in H^\rho_0(\Omega;\R^n)$ such that 
\begin{equation*}\int_{\R^n} C[D_{\rho,sym}v]:D_{\rho,sym}w\,dx=\int_{\R^n} f\cdot w\,dx,\,\forall w\in H^{\rho}_0(\Omega;\R^n).\end{equation*}

Of particular interest is the isotropic case, when $W$ is isotropic, the elasticity tensor $C$ takes the form (see \cite{Knowles95} and the references therein) $$C[M]=2\mu M+\lambda(\operatorname{tr}M)I,\,M\in\R^{n\times n}_{sym},$$ where $\mu,\lambda$ are the Lam\'e's moduli of the material (which represents the bulk and shear modulus respectively in the context of solid materials). In order to the tensor $C$ to satisfy the strong ellipticity condition, we require to $$M:C[M]>0,\quad \forall M=a\otimes b,\,a,b\in \R^n\setminus\{0\}.$$ 
Recall that for vectors $a,b\in \R^n$, the tensor product $a\otimes b$ is the matrix $(a_ib_j)_{i,j}$.\\
For such $M$, \begin{align*}M:C[M]&=M:C[M_{sym}]=(a\otimes b):\frac{1}{2} C[a\otimes b+b\otimes a]\\
&=\mu\left((a\otimes b):(a\otimes b)+(a\otimes b):(b\otimes a)\right)+\lambda \operatorname{tr}(a\otimes b)(a\otimes b):I.
\end{align*}
Since \begin{align*}
    (a\otimes b):(a\otimes b)&=(a\cdot a)(b\cdot b)=|a|^2|b|^2,\\
    \operatorname{tr(a\otimes b)}&=a\cdot b,\\
    (a\otimes b):I&=\operatorname{tr}(a\otimes b)=a\cdot b,\\
    (a\otimes b):(b\otimes a)&=(a\cdot b)(b\cdot a)=(a\cdot b)^2,
\end{align*}
we have \begin{align*}M:C[M]&=\mu\left(|a|^2|b|^2+(a\cdot b)^2\right)+\lambda (a\cdot b)^2\\
&=\mu\left(|a|^2|b|^2-(a\cdot b)^2\right)+(\lambda+2\mu) (a\cdot b)^2\end{align*} Now, for $\mu>0$ and $2\mu+\lambda>0$, implies that $C$ is strongly elliptic. Conversely, if $C$ is strongly elliptic, taking $a=b$ yields that $\lambda+2\mu>0$ and taking $a$ orthogonal to $b$ implies that $\mu>0$. Hence, $$C\,\text{is strongly elliptic}\iff \begin{cases}
    \mu>0,\\
    2\mu+\lambda>0
\end{cases}.$$

To write the equations for this particular case we require the following technical lemma, which relates the nonlocal gradient with the nonlocal divergence as it happens in the classical setting.
\begin{Lema}\label{Traza}
    Let $u\in H^{\rho,p}_0(\Omega;\R^n)$. Then, $$\operatorname{tr}D_{\rho,sym}u=\operatorname{tr}(D_\rho u)^T=\operatorname{tr}D_\rho u=\operatorname{div}_\rho u$$
\end{Lema}
The proof of this lemma is analogous to \cite[Lemma 2.5]{BellidoCuetoMoraCorral2021}.

In this case, $$C[D_{\rho,sym}v]=2\mu D_{\rho,sym}v+\lambda(\operatorname{tr}D_{\rho,sym}v)I=2\mu D_{\rho,sym}v+\lambda(\operatorname{div}_\rho v)I,$$ by the previous lemma. Now, $$C[D_{\rho,sym}v]:D_{\rho,sym}w=2\mu D_{\rho,sym}v:D_{\rho,symw}+\lambda(\operatorname{div}_\rho v)I:D_{\rho,sym}w,$$ 
where it is straightforward to see that $$I:D_{\rho,sym}w=\operatorname{tr}D_{\rho,sym}v=\operatorname{div}_\rho w,$$ and hence the weak form of the isotropic nonlocal linear elasticity for general kernels is: find $v\in H_0^\rho(\Omega;\R^n)$ such that 
$$\int_{\R^n}\left(2\mu D_{\rho,sym}v:D_{\rho,sym}w+\lambda\operatorname{div}_\rho v\operatorname{div}_\rho w\right)\,dx=\int_{\R^n}f\cdot w,\,\forall w\in H_0^\rho(\Omega;\R^n).$$ In order to write the strong form we need further results regarding calculus for nonlocal operators with general kernels. First of all we recall the definition of the following Laplacian-type operator introduced in \cite{GarciaSaez2026} in analogy to the local and fractional cases:
\begin{Defi}\em
    Let $u\in C_c^\infty(\R^n)$ or $u\in C_c^\infty(\R^{n\times n})$. We define the $\rho$-nonlocal Laplacian as $$(-\Delta)_\rho u:=-\operatorname{div}_\rho D_\rho u.$$
\end{Defi} 
Let us assume that $v\in C_c^\infty(\Omega;\R^n)$. We have that $$\operatorname{div}_\rho(C[D_{\rho,sym}v])=2\mu \operatorname{div}_\rho (D_{\rho,sym}v)+\lambda \operatorname{div}_\rho\left((\operatorname{div}_\rho v)I\right).$$ In \cite[Lemma 3.8]{BellidoMoraHidde2024} a nonlocal Leibniz rule was proven for $D_\rho$. Briefly adapting the proof we can get the following result.
\begin{Lema}\label{NLLeib}
    Let $1<p<\infty$, $\varphi\in C_c^\infty(\R^n)$ and $\Phi\in H^{\rho,p}(\R^n;\R^{n\times n})$. Then, $$\operatorname{div}_\rho (\varphi\Phi)=\varphi\operatorname{div}_\rho \Phi+K_{\rho,\varphi}(\Phi),$$ where $$K_{\rho,\varphi}(\Phi)(x):=\int_{\R^n}\frac{\varphi(x)-\varphi(y)}{|x-y|}\Phi(y)\frac{x-y}{|x-y|}\rho(x-y)\,dy.$$
    \end{Lema}
In this case, taking $\varphi=\operatorname{div}_\rho v$ y $\Phi=I$, yields \begin{align*}\operatorname{div}_\rho\left((\operatorname{div}_\rho v)I\right)&=\operatorname{div}_\rho v\operatorname{div}_\rho I+K_{\rho,\operatorname{div}_\rho}(I)=K_{\rho,\operatorname{div}_\rho v}(I)\\
&=\int_{\R^n}\frac{\operatorname{div}_\rho v(x)-\operatorname{div}_\rho v(y)}{|x-y|}\frac{x-y}{|x-y|}\rho(x-y)\,dy\\
&=D_\rho(\operatorname{div}_\rho v).\end{align*}
On the other hand, $$\operatorname{div}_\rho (D_{\rho,sym}v)=\frac{\operatorname{div}_\rho (D_\rho v)+\operatorname{div}_\rho (D_\rho v)^T}{2}=\frac{-(-\Delta)_\rho v+\operatorname{div}_\rho (D_\rho v)^T}{2}.$$

Now, noticing that $\operatorname{div}_\rho (D_\rho v)^T=D_\rho (\operatorname{div}_\rho v)$, for all $u\in C_c^\infty(\R^n;\R^n)$, result that generalizes \cite[Lemma 3.4]{BellidoCuetoMora2023} and can be proved in an analogous way, we find that $$\operatorname{div}_\rho (C[D_{\rho,sym}v])=-\mu(-\Delta)_\rho v+(\mu+\lambda)D_\rho(\operatorname{div}_\rho v),$$ and, hence, the strong form of the nonlocal linear equations for general kernels in the isotropic case reads as: $$\begin{cases}
    \mu(-\Delta)_\rho v-(\mu+\lambda)D_\rho(\operatorname{div}_\rho v)&=f,\,\text{in $\Omega$},\\
    \hfill v&=0,\,\text{on $\Omega^c$}.
\end{cases}$$

Choosing $\rho(x)=\rho^s(x)=c_{n,s}\frac{1}{|x|^{n+s-1}}, s\in (0,1)$, we arrive at the fractional isotropic problem $$\begin{cases}
    \mu(-\Delta)^s v-(\mu+\lambda)D^s(\operatorname{div}^s v)&=f,\,\text{in $\Omega$},\\
    \hfill v&=0,\,\text{on $\Omega^c$}.
\end{cases},$$ where $(-\Delta)^s$ is the widely-used fractional Laplacian (see \cite{Kwasnicki2017}, for a survey on this operator). We point the strong relationship of our model with he fractional Lam\'e-Navier model studied in \cite{Scott2022}. In that work, the fractional model is obtained as the $s$-fractional power the classical Lam\'e-Navier system arriving at 
$$\begin{cases}
    \mu^s(-\Delta)^s v-\left((2\mu+\lambda)^s-\mu^s\right)D^s(\operatorname{div}^s v)&=f,\,\text{in $\Omega$},\\
    \hfill v&=0,\,\text{on $\Omega^c$}.
\end{cases}.$$ 
Notice that this model coincides with our model once in the fractional case $\rho=\rho^s$, once we find the right assingment for the parameters in both models.

\section{Existence and uniqueness of solutions}

In this, we prove well-posedness of the linearized model obtained in the previous section. We consider both the Dirichlet and Neumann boundary conditions. Although, the analysis for Neumann boundary conditions for nonlocal problems in its infancy, we can consider the particular case in which $\rho=\rho_\delta^s$, following \cite{Kreisbeck2024}. 

\subsection{Dirichlet boundary conditions for general nonlocal gradients}

The two key ingredients for the existence and uniqueness of weak solutions for the problem (\ref{eqn:P}) are a nonlocal version for Korn's inequality for general kernels and the Poincar\'e inequality (Theorem \ref{PoincareII}). In order to prove the former, we will translate the classical Korn's inequality using the nonlocal to local operator from Lemma \ref{NLTOL}. 
\begin{Lema}\label{Korn}
    Let $1<p<\infty$. Then, there exists a constant $c_p>0$ only depending on $p$ such that for all $u\in W^{1,p}(\R^n;\R^n)$, $$\norm{D_{sym}u}_{L^p(\R^n)}\geq c_p\norm{Du}_{L^p(\R^n)}.$$
\end{Lema}
See \cite[Lemma 7.1]{BellidoCuetoMora2023} for a proof. Now, we can prove the general nonlocal Korn's inequality.
\begin{Prop}\label{NLKorn}
    Let $1<p<\infty$. Then, there exists $c_p>0$ such that for all $u\in H^{\rho,p}(\R^n;\R^n)$, $$\norm{D_{\rho,sym}u}_{L^p(\R^n)}\geq c_p\norm{D_\rho u}_{L^p(\R^n)}.$$
\end{Prop}
\noindent{}\textbf{Proof:} By density, it is enough to prove the result for functions $u\in C_c^\infty(\R^n;\R^n)$. For such $u$, we define the function $v=\mathcal{Q}_\rho u=Q_\rho*u$. Then, by Lemma \ref{NLTOL}, $v\in W^{1,p}(\R^n;\R^n)$ and $Dv=D_\rho u$. Then, by Lemma \ref{Korn}, $$\norm{D_{\rho,sym}u}_{L^p(\R^n)}=\norm{D_{sym}v}_{L^p(\R^n)}\geq c_p\norm{Dv}_{L^p(\R^n)}=c_p\norm{D_{\rho}u}_{L^p(\R^n)},$$ as we wanted to prove. \qed

Proved this, we are ready to establish de existence and uniqueness of weak solutions for the problem (\ref{eqn:P}).
\begin{Teor}\label{Existence}
    Assume $C:\R^{n\times n}_{sym}\to\R^{n\times n}_{sym}$ is positive definite. Let $f\in L^2(\Omega;\R^n)$. Then, there exists a unique weak solution to  (\ref{eqn:P}). Moreover, it is the unique minimizer of the functional $$E(v):=\frac{1}{2}\int_{\R^n}C[D_{\rho,sym}v]:D_{\rho,sym}v\,dx-\int_{\Omega}f\cdot v,\,v\in H^{\rho}_0(\Omega;\R^n).$$
\end{Teor}
\noindent{}\textbf{Proof:} This result is just a standard application of the classical Lax-Milgram's lemma. Clearly, the bilinear form $$a(v,w):= \int_{\R^n}C[D_{\rho,sym}v]:D_{\rho,sym}w\,dx,$$ is continuous in $H^{\rho}_0(\Omega;\R^n)$. Moreover, by the hypothesis under $C$, it is symmetric. The linear form $$w\mapsto \int_\Omega f\cdot w,$$ is also continuous in $H^{\rho}_0(\Omega;\R^n)$, hence in order to prove the result it is enough to check the coercivity of the bilinear form. Note that since $C$ is positive definite in $\R^{n\times n}_{sym}$, there exists $c_1>0$ such that for every $M\in \R^{n\times n}_{sym}$, $$C[M]:M\geq c_1|M|^2.$$ Hence, for every $v\in H_0^\rho(\Omega;\R^n)$, by Lemma \ref{NLKorn} and Theorem \ref{PoincareII}, \begin{align*}
a(v,v)&=\int_{\R^n}C[D_{\rho,sym}v]:D_{\rho,sym}v\,dx\geq c_1\int_{\R^n}|D_{\rho,sym}v|^2\,dx=c_1\norm{D_{\rho,sym}v}_{L^2(\R^n)}^2\\
&\geq c_2\norm{D_\rho v}_{L^2(\R^n)}^2\geq c_2'\norm{u}_{H^{\rho}_0(\Omega;\R^n)},
\end{align*} for some positive constant $c_2'$. Hence, the bilinear form $a$ is coercive and the result follows. \qed

In fact, it is enough for $f$ to lie in the dual of $H^{\rho}_0(\Omega;\R^n)$, which is a space that has not been studiet yet, but we know that it contains $L^q(\Omega;\R^n)$ for $$\begin{cases}
    q\in [\frac{2n}{2s+n},\infty],\,2t<n,\\
    q\in [1,\infty],\,2s>n.
\end{cases},$$ thanks to the continuous embeddings for $H^{\rho,p}_0(\Omega;\R^n)$, with $s,t$ being as in the hypothesis $(H0)-(H4)$.

In the fractional case, $\rho^s=c_{n,s}\frac{1}{|\cdot|^{n-s+1}}$, yielding to Riesz fractional gradient, $D_{\rho^s}=D^s$ and $H^{\rho^s}_0(\Omega;\R^n)=H^s_0(\Omega;\R^n)$, the model corresponds to the linear fractional elasticity equations studied in \cite[Section 7]{BellidoCuetoMora2023} and \cite{Silhavy2024}. Moreover, when the kernel is taken as $\rho^s_\delta(x)=w_\delta(x)\frac{1}{|x|^{n-s+1}}$, for some $\delta>0$ small enough and $w_\delta$ a proper cut-off function we obtain the truncated Riesz fractional gradient, $D_\rho=D^s_\delta$, the space $H^{\rho_\delta^s}_0(\Omega;\R^n)=H^{s,2,\delta}_0(\Omega;\R^n)$, introduced in \cite{BellidoCuetoMoraCorral2023}, which can be characterized as the space $$\{u\in H^{s,2,\delta}(\Omega;\R^n): u=0\,\text{a.e. in}\,\Omega_\delta\setminus\Omega_{-\delta}\},$$ as it was proven in \cite[Proposition 3]{Cueto2023}, where $$H^{s,2,\delta}(\Omega;\R^n)=\{u\in L^2(\Omega_\delta;\R^n): D^{s}_\delta u\in L^2(\Omega;\R^{n\times n})\}.$$  
In this case the external force $f$ is taken in $L^2(\Omega;\R^n)$ since at that time the dual space of $H^{s,\delta}_0(\Omega;\R^n)$ was not known. However, by \cite[Lemma 5]{Cueto2023}, we known that whenever $\Omega_{-\delta}$ has Lipschitz boundary, $$H^{s,p,\delta}_0(\Omega;\R^n)=H^{s,p}_0(\Omega_{-\delta};\R^n),\,1<p<\infty,$$ with equivalence of norms. Hence, it is enough to take $f$ in the space $$\left(H^{s,p,\delta}_0(\Omega;\R^n)\right)^*=\left(H^{s,p}_0(\Omega_{-\delta};\R^n)\right)^*=H^{-s,p'}(\Omega_{-\delta};\R^n),$$ with $1/p+1/p'=1$.\\

In view of the previous observations, the existence and uniqueness result of Theorem \ref{Existence} for the problem \ref{eqn:P}, generalizes \cite[Theorem 6.2, Theorem 7.3]{BellidoCuetoMora2023}.

\subsection{Neumann boundary conditions for $\rho=\rho^s_\delta$}

As we have mentioned above, a general theory for Neumann boundary conditions for problems involving general nonlocal gradients is still pending. However, in the particular case $\rho=\rho_\delta^s$ this has been developed in \cite{Kreisbeck2024}. We follow this reference along this section.   

Let fix $\Omega$ a bounded Lipschitz domain, $1<p<\infty$ and we consider the energy integrand $I_\delta^s:H^{s,p,\delta}(\Omega;\R^n)\to \R\cup\{\infty\}$, defined as $$I_\delta^s[u]:=\int_{\Omega}W(x,D^s_\delta u)\,dx-\int_{\Omega_\delta}f\cdot u\,dx,$$ where $W:\Omega\times \R^{n\times n}\to \R\cup\{\infty\}$ is a suitable Carath\'eodory integrand and $f\in L^{p'}(\Omega_\delta;\R^n)$. We suppose that $W$ is under the hypothesis of \cite[Theorem 6.1]{Kreisbeck2024} which establishes the existence of minimizers for $I_\delta^s$ in the subspace $N^{s,p,\delta}(\Omega)^\perp$, and allows us
to differentiate under the integral sign in order to obtain the Euler-Lagrange equations. We further assume the compatibility condition $$\int_{\Omega_\delta}f\cdot v\,dx=0,\,v\in {N}^{s,p,\delta}(\Omega).$$ Now, we can obtain the Euler-Lagrange equations for the model as we did in Section 3. The difference will be that since we are not imposing that the value of the functions are zero outside $\Omega$, nonlocal Neumann conditions will naturally appear. We have that 
$$\frac{d}{dt}I_\delta^s[u+tv]\bigg|_{t=0}=\int_\Omega W_2(x,D_\delta^su)\cdot D_\delta^sv\,dx-\int_{\Omega_\delta}f\cdot v\,dx,$$ 
and by the nonlocal integration by parts formula \cite[Theorem 3.11]{BellidoCueto2024},
$$-\int_{\Omega_{-\delta}}\operatorname{div}_\delta^s\left(W_2(x,D_\delta^s u)\right)v\,dx+\int_{\Gamma_{\pm\delta}}\mathcal{N}_\delta^s\left(W_2(x,D_\delta^s u)\right)\cdot v=\int_{\Omega_\delta}f\cdot v\,dx.$$ 
Finally, applying the fundamental lemma of the calculus of variations we obtain the nonlocal Neumann problem 
$$\begin{cases}
    -\operatorname{div}_\delta^s W_2(x,D_\delta^s u)&=f,\,\text{in $\Omega_{-\delta}$},\\
    \mathcal{N}_\delta^s \left(W_2(x,D_\delta^s u)\right)&=f,\,\text{on $\Gamma_{\pm \delta}$},
\end{cases}$$
(see \cite{Kreisbeck2024}). By simplicity, we may assume that $f|_{\Gamma_{\pm \delta}}=0$. We can perform the same linearization process as above, arriving at  the strong form of the equations of general nonlocal linear elasticity with nonlocal Neumann conditions,
\begin{equation}\begin{cases}\label{eqn:PN}
    -\operatorname{div}_\delta^s(C[D_{\delta,sym}^su])&=f,\,\text{in $\Omega_{-\delta}$},\\
    \mathcal{N}_\delta^s(C[D_{\delta,sym}^su])&=0,\,\text{on $\Gamma_{\pm\delta}$}.
\end{cases}\end{equation}
The problem in its weak formulation reads as: find $u\in N^{s,2,\delta}(\Omega;\R^n)^\perp$, such that $$\int_\Omega CD_{\delta,sym}^su:D^s_{\delta,sym}v=\int_{\Omega _\delta}f \cdot v,\quad \forall v\in N^{s,2,\delta}(\Omega;\R^n)^\perp.$$
We are ready to establish the existence and uniqueness of weak solutions for this problem.
\begin{Teor}\label{ExistenceNeumann}
    Assume $C:\R^{n\times n}_{sym}\to\R^{n\times n}_{sym}$ is positive definite, and $f\in L^2(\Omega_\delta;\R^n)$. Then, there exists a unique weak solution to  (\ref{eqn:PN}). Moreover, it is the unique minimizer of the functional $$E(u):=\frac{1}{2}\int_{\Omega}CD_{\delta,sym}^su:D^s_{\delta,sym}u-\int_{\Omega_\delta}f\cdot u,\,u\in N^{s,2,\delta}(\Omega;\R^n)^\perp.$$
    
\end{Teor}
\noindent{}\textbf{Proof:} Again, the result follows as a standard application of the classical Lax-Milgram's lemma. Clearly, the bilinear form $$a(v,w):= \int_{\R^n}C[D^s_{\delta,sym}v]:D^s_{\delta,sym}w\,dx,$$ is continuous in $N^{s,2,\delta}(\Omega;\R^n)^\perp$. Moreover, by the hypothesis under $C$, it is symmetric. The linear form $$w\mapsto \int_\Omega f\cdot w,$$ is also continuous in $N^{s,2,\delta}(\Omega;\R^n)^\perp$, hence in order to prove the result it is enough to check the coercivity of the bilinear form. Since $C$ is positive definite in $\R^{n\times n}_{sym}$, there exists $c_1>0$ such that for every $M\in \R^{n\times n}_{sym}$, $$C[M]:M\geq c_1|M|^2.$$ Hence, for every $v\in N^{s,p,\delta}(\Omega;\R^n)^{\perp}$, by Lemma \ref{NLKorn} (see also \cite[Proposition 6.1]{BellidoCuetoMora2023}) and Lemma \ref{PW}, \begin{align*}
a(v,v)&=\int_{\Omega}C[D^s_{\delta,sym}v]:D^s_{\delta,sym}v\,dx\geq c_1\int_{\Omega}|D^s_{\delta,sym}v|^2\,dx=c_1\norm{D^s_{\delta,sym}v}_{L^2(\Omega;\R^n)}^2\\
&\geq c_2\norm{D_\delta^s v}_{L^2(\Omega;\R^n)}^2\geq c_2'\norm{v-\pi_\delta^s(v)}_{L^{2}(\Omega_\delta;\R^n)}^2=c_2'\norm{v}_{L^{2}(\Omega_\delta;\R^n)}^2,
\end{align*} for some positive constant $c_2'$. Hence, $$a(v,v)\geq C\norm{v}_{H^{s,2,\delta}(\Omega;\R^n)}^2,$$ so the bilinear form $a$ is coercive and the result follows. \qed

Notice that if the data $f$ satisfies
$$\int_{\Omega_\delta}f\cdot h=0,\quad \forall h\in N^{s,2,\delta}(\Omega;\R^n),$$
then the energy is invariant under translations in $N^{s,2,\delta}(\Omega;\mathbb{R}^n)$, and consequently, existence of solutions in the quotient $H^{s,2,\delta}(\Omega;\mathbb{R}^n)/N^{s,2,\delta}(\Omega;\mathbb{R}^n)$ holds.

\section{Connection with Eringen's model}

One of the most popular models, particularly due to its application in engineering, that employs a nonlocal perspective for continuum mechanics is the Eringen model \cite{Eringen2002}. In the particular case of linear elasticity, nonlocality enters into the model through stress filtering (see \cite{Polizzotto2001}); namely, if we call $v$ the displacement, then stress is computed from $v$ as
\[\sigma (x) =\int_{\Omega}A(x,y)C[D_{sym} v(y)] \, dy, \text{  in $\Omega$},\]
where $A$ is the nonlocal interaction kernel between particles accounting for stress interactions, and the equilibrium equation is
\[ -\operatorname{div} \sigma =f,  \text{  in $\Omega$}.\]
Nonlocality kernel $A$ is typically assumed to be positive definite, i.e., it satisfies the so-called Mercer's condition, $$\int_{\Omega\times \Omega}A(x,y)\phi(x)\phi(y)\,dx\,dy>0,$$
for any test function $\phi$. We additionally assume that there exists a function $\tilde{A}:(0,\infty)\to [0,\infty)$, such that $A(x,y)=\tilde{A}(|x-y|)$. Here, the function $\tilde{A}(d)$ is seen as a description of nonlocal interaction of particles at a distance $d>0$, and we are assuming radial interaction indeed.

In \cite{BellidoCuetoMora2023}, the relationship among linear elasticity models based on nonlocal gradients and the Eringen model was explored, finding the somewhat surprising fact that there are particular elections of $\tilde A$, concretely
$$\tilde{A}(|x|)=I_{2-2s}(x)=I_{1-s}*I_{1-s}(x), \text{  and  }\tilde{A}(|x|)=Q_\delta^s*Q_\delta^s(x),$$
making he homogeneous Eringen model equivalent to \eqref{eqn:P} for $\rho=\rho^s$ and $\rho=\rho_\delta^s$ respectively, i.e. models based on the Riesz-fractional gradient and the truncated Riesz-fractional gradient respectively. This relationship may be generalize for general nonlocal gradients in a more or less straightforward manner following \cite{BellidoCuetoMora2023} that we include in this section for the sake of completeness.

For simplicity, we assume homogeneous Dirichlet boundary conditions. Following \cite{Polizzotto2001}, nonlocal Eringen model is
 \begin{equation} \label{eq: Eringen model}
  \begin{cases}
  -\operatorname{div} \sigma =f,  &\text{ in $\Omega$} \\
  \sigma (x) =\int_{\Omega}A(x,y)C[D_{sym} v(y)] \, dy, &\text{ in $\Omega$} \\
  v=0, &\text{ on $\partial\Omega $} .
  \end{cases}  
  \end{equation}
Here $C$ is the constant elasticity tensor. Let us define the symmetric bilinear form $$(v,w)_A:=\int_{\Omega\times \Omega}A(x,y)Dv(x):Dv(y)\,dx\,dy,$$ which defines an inner product over $C_c^\infty(\R^n)$, and hence, induces a norm, which we will denote as $\norm{\cdot}_A$. We define the space $$V_A(\Omega;\R^n):=\overline{C_c^\infty(\Omega;\R^n)}^{\norm{\cdot}_A},$$
where with a small abuse of notation, we assume that functions are extended to $\R^n$ by zero outside $\Omega$. Obtaining the weak formulation of (\ref{eq: Eringen model}) follows as usual, arriving at the problem: find $v\in V_A(\Omega;\R^n)$ such that $$\int_{\Omega\times \Omega}A(x,y)C[D_{sym}v](x):D_{sym}w(y)\,dy\,dx=\int_\Omega f\cdot w,\,\forall w\in V_{A}(\Omega;\R^n).$$

In this section we assume that 
$$\tilde{A}(|x|)=Q_\rho (x)*Q_\rho(x),$$ 
for a general kernel $\rho$, and we will show that this choice makes \eqref{eq: Eringen model} equivalent to \eqref{eqn:P}, generalizing the fractional gradient cases. In order to prove that, we have to establish the equivalence of the spaces $V_A(\Omega;\R^n)$ and $H^\rho_0(\Omega;\R^n)$, and the equality of the bilinear form $$a(v,w)_A:=\int_{\R^n}\int_{\R^n}\tilde{A}(|x-y|)C[D_{sym}v](x):D_{sym}w(y)\,dy\,dx,$$ with $$ \int_{\R^n}C[D_{\rho,sym}v](x):D_{\rho,sym}w(x)\,dx.$$
Note that this choice for $\tilde{A}$ is valid since $Q_\rho$ is radial and the convolution of radial functions is radial. Hence, $A_\rho(x,y)=\tilde{A}_\rho(|x-y|)$ depends on the modulus $|x-y|$ and is strictly positive, as well as it lies in $L^1(\R^n\times \R^n)$ by Young's convolution inequality.
To emphasize the choice of $\tilde{A}(|x|)$ as $Q_\rho*Q_\rho(x)$, we rename the space as $V_{A_\rho}(\Omega;\R^n)$ and the bilinear form as $a(v,w)_{A_\rho}$. The main result of the section is the following: 

\begin{Teor}\label{Eringen}
Let $A_\rho(x,y)$ defined as $A_\rho(x,y):=\tilde{A}_\rho(|x-y|)$, where $\tilde{A}_\rho(|x|)=Q_\rho*Q_\rho(x)$. Then, $$V_{A_\rho}(\Omega;\R^n)=H^\rho_0(\Omega;\R^n),$$ with equivalence of norms, and $$a(v,w)_{A_\rho}=\int_{\R^n}C[D_{\rho,sym}v](x):D_{\rho,sym}w(x)\,dx.$$ Hence, the problems \eqref{eqn:P} and \eqref{eq: Eringen model} are equivalent in their weak formulation.
\end{Teor}
\noindent{}\textbf{Proof:}
By density, it is enough to prove that $(\cdot,\cdot)_{A_\rho}$ induces an equivalent norm to the $H^\rho$-norm for functions $C_c^\infty(\Omega;\R^n)$. Let $v\in C_c^\infty(\Omega;\R^n)$. Then, \begin{align*}
    (v,v)_{A_\rho}&=\int_{\R^n}\int_{\R^n}A_\rho(x,y)Dv(x):Dw(y)\,dx\,dy\\
    &=\int_{\R^n}\int_{\R^n}(Q_\rho*Q_\rho)(x-y)Dv(x):Dv(y)\,dx\,dy\\
    &=\int_{\R^n}\left[(Q_\rho*Q_\rho)*Dv\right](y):Dv(y)\,dy,
\end{align*} and by \cite[Proposition 4.16]{Brezis2010}, $$\int_{\R^n}\left[(Q_\rho*Q_\rho)*Dv\right](y):Dv(y)\,dy=\int_{\R^n}(Q_\rho*Dv)(y):(Q_\rho*Dv)(y)\,dy=\int_{\R^n}|Q_\rho*Dv(y)|^2\,dy,$$ 
and since by \cite[Proposition 2.6]{BellidoMoraHidde2024}, $D_\rho v=Q_\rho*Dv$, $$(v,v)_{A_\rho}=\norm{D_\rho v}^2_{L^2(\R^n;\R^{n\times n})}.$$ Now, by Theorem \ref{Poincare}, $\norm{D_\rho v}_{L^2(\R^n;\R^{n\times n})}$ is equivalent to the $H^\rho$-norm on the subspace $H^\rho_0(\Omega;\R^n)$, and hence, this proves the equality  $$V_{A_\rho}(\Omega;\R^n)=H^\rho_0(\Omega;\R^n),$$ with equivalence of norms. By an analogous argument, \begin{align*}
    a(v,w)_{A_\rho}&=\int_{\R^n}\int_{\R^n}\tilde{A}(|x-y|)C[D_{sym}v](x):D_{sym}w(y)\,dx\,dy\\
    &=\int_{\R^n}\int_{\R^n}(Q_\rho*Q_\rho)(x-y)C[D_{sym}v](x):D_{sym}w(y)\,dx\,dy\\
    &=\int_{\R^n}\left[(Q_\rho*Q_\rho)*C[D_{sym}v]\right](y):D_{sym}w(y)\,dy\\
    &=\int_{\R^n}(Q_\rho*C[D_{sym}v])(y):(Q_\rho *D_{sym}w)(y)\,dy\\
&=\int_{\R^n}C[D_{\rho,sym}v](y):D_{\rho,sym}w(y)\,dy,
\end{align*} as we wanted to prove.\qed

Notice that, since the existence and uniqueness of solutions for the problem \eqref{eqn:P} hold, there exist weak solutions for the problem \eqref{eq: Eringen model}. This is interesting, as there are no results on the existence and uniqueness of solutions for the Eringen model in the literature, as pointed out in \cite{EvgrafovBellido2019}.



 

\section{Localization results}

This section is devoted to analyzing the convergence of nonlocal models to their local, or fractional, counterparts. We distinguish between the Dirichlet and Neumann cases.

\subsection{Localization for the Dirichet problem}

We proceed in two directions. First, we study homogenization with respect to the kernel $\rho$; that is, we consider problems associated with the rescaled kernel $\rho_\delta$ and analyze the limits as the parameter $\delta$ tends to zero and to infinity. Second, we address the case $\rho = \rho_\delta^s$, which gives rise to the truncated Riesz fractional gradient, and we investigate the limit as $s \to 0^+$.

The material in this section is not new. Our analysis relies entirely on the results in \cite{CuetoHidde2025} concerning homogenization for general gradients, and on \cite{Cueto2023} regarding localization in the fractional parameter. Our goal here is not to claim novelty, but rather to apply these existing results within the framework of nonlocal linear elasticity, thereby contributing to the completeness of the theory.

\subsubsection{Homogenization for general gradients}

We consider an arbitrary kernel $\rho$ verifying $(H0)$--$(H4)$, such that $\operatorname{supp}\rho=\overline{B(0,1)}$. Given a horizon $\delta>0$ such that $\Omega_{-\delta}$ is well defined, the rescaled kernel $\rho_\delta$ is defined as
\[
\rho_\delta(x):=c_\delta\rho(x/\delta), \quad x\in \R^n\setminus\{0\}.
\]
The choice of the scaling factor implies that
\[
\widehat{Q_{\rho_\delta}}(\xi)=\widehat{Q_\rho}(\delta\xi).
\]
We can address the localization problems  in $H^{\rho,p,\delta}_0(\Omega)$ as $\delta\to 0^+$ and $\delta\to +\infty$, and the constant $c_\delta$ will specify later depending on asymptotics for $\delta$. With the observation that the Dirichlet condition in $H^{\rho,p,\delta}_0(\Omega)$ is equivalent to have $u=0$ in $\Gamma_{\pm \delta}$, we will rename the domain in order to make meaningful both limit passages $\delta\to 0^+$ and $\delta\to+\infty$. In fact, we will prescribe Dirichlet conditions on $\Omega_{2\delta}\setminus\Omega$, which is clearly equivalent to the previous one. In this case, problem \eqref{eqn:P} reads as \begin{equation}\label{eqn:Pdelta}
    \begin{cases}
    -\operatorname{div}_{\rho_\delta} \left(C[D_{\rho_\delta,sym}v]\right)&=f,\,\text{in $\Omega$},\\
    \hfill v&=0,\,\text{on $\Omega_{2\delta}\setminus\Omega$}.
\end{cases}
\end{equation}

For the case of vanishing horizon, i.e., the limit $\delta\to 0^+$, we consider the scaling factors $c_\delta$ to be $c_\delta=\delta^{-n},$ since they preserve the normalizations $$\int_{\R^n}\rho_\delta\,dx=n,\,\int_{\R^n}Q_{\rho_\delta}\,dx=1.$$
Let $W:\Omega\times \R^{n\times n}\to \R$ a suitable Carath\'eodory integrand, and associated to it the family of vectorial energy functionals $(I_\delta)_{\delta\in (0,1]}$, with $I_\delta: L^p(\R^n;\R^n)\to \R\cup\{\infty\}$, defined as $$I_\delta(u):=\begin{cases}
 \displaystyle   \int_{\Omega_{\delta}}W(x,D_{\rho_\delta}u)\,dx-\int_{\Omega_\delta}f\cdot u\,dx, & u\in H_0^{\rho,p,\delta}(\Omega;\R^n)\\
    \infty,&\text{otherwise},
\end{cases}$$ for $f\in L^{p'}(\Omega_\delta;\R^n)$, where we consider that the functions in the domain of $I_\delta$ are defined on $\R^n$ with zero Dirichlet conditions on $\Omega^c$. Assuming the following conditions on the integrand $W$,
\begin{itemize}
    \item $c_1|A|^p-c_2\leq W(x,A)\leq c_2(1+|A|^p),\,x\in \Omega,\,A\in \R^{n\times n},$
    \item $W(x,\cdot)$ is quasiconex a.e. $x\in \Omega$,
\end{itemize} we have the following result \cite[Theorem 3.7]{CuetoHidde2025}.

\begin{Teor}
    Let $$I_0(u):=\begin{cases}
    \displaystyle\int_{\Omega}W(x,Du)\,dx-\int_{\Omega}f\cdot u\,dx, & u\in W_0^{1,p}(\Omega;\R^n),\\
    \infty,& \text{otherwise}.
\end{cases}$$ Then, the family $(I_\delta)_\delta$, $\Gamma$-converges with respect to the $L^p(\R^n;\R^n)$-convergence to the functional $I_0$ as $\delta\to 0^+$, $$\Gamma(L^p)-\lim_{\delta\to 0^+}I_\delta=I_0.$$ Additionally, $(I_\delta)_\delta$ is equi-coercive with respect to strong convergence in $L^p(\R^n;\R^n)$.
\end{Teor}
A direct consequence of this result is the $\Gamma$-convergence of the nonlocal problem with finite horizon \eqref{eqn:Pdelta} to the classical local Lamé problem of linearized elasticity. In fact, we have that $$\Gamma(L^2)-\lim_{\delta\to 0^+}\begin{cases}
    -\operatorname{div}_{\rho_\delta} \left(C[D_{\rho_\delta,sym}]v\right)&=f,\,\text{in $\Omega$},\\
    \hfill v&=0,\,\text{on $\Omega_{2\delta}\setminus\Omega$}.
\end{cases}=\begin{cases}
    -\operatorname{div}\left(C[D_{sym}v]\right)&=f,\,\text{in $\Omega$},\\
    \hfill v&=0,\,\text{on $\partial\Omega$}.
\end{cases}.$$ The equi-coercivity of the family $(I_\delta)_\delta$ and the existence and uniqueness of weak solutions for the nonlocal problem in $H^{\rho,2,\delta}_0(\Omega;\R^n)$ implies the existence and uniqueness of weak solutions for the local one in $W^{1,2}_0(\Omega;\R^n)$. 

On the other hand, it is natural to expect that when the horizon of interaction diverges, and hence the domain $\Omega_\delta$ becomes the whole $\R^n$,  the gradient $D_\delta^s$ converges to a purely fractional one $D^{s_\infty}$. Surprisingly, this is also the case not just for $D_\delta^s$, but for every rescaled family $\rho
_\delta$ of general kernels under our hypothesis, as it was proven in \cite{CuetoHidde2025}.\\

In order to overcome the lack of integrability of the fractional kernel in $\R^n$, the rescaling constants $c_\delta$ are chosen as $$c_\delta=\overline{\rho}(1/\delta)^{-1},$$
where $\bar\rho$ is the radial representative of kernel $\rho$. We also require the pointwise convergence of $\rho_\delta$ on $\R^n\setminus\{0\}$ when $\delta\to +\infty$, to some $\rho_\infty$. The key is that this pointwise limit is indeed a fractional kernel \cite[Lemma 4.2]{CuetoHidde2025}, i.e., there exists some $s_\infty\in [s,t]$, such that $$\rho_\infty(x)=\frac{1}{|x|^{n+s_\infty-1}},\,x\not =0,$$ where $s,t$ are the fractional parameters of the hypothesis $(H3)-(H4)$. In fact, as it is mentioned in \cite[Remark 4.3]{CuetoHidde2025}, the exact value of $s_\infty$ can be computed as $$s_\infty=\log\left(\overline{\rho}_\infty(1/e)\right)-n+1.$$
Consider now a Carath\'eodory integrand $W:\R^n\times \R^{n\times n}\to \R$ and the family of energy functional $(I_\delta)_\delta$ defined as in the previous section. We also consider the functional $I_\infty$, defined as $$I_\infty(u):=\begin{cases}
   \displaystyle \int_{\R^n}W(x,D^{s_\infty}u)\,dx-\int_{\R^n}f\cdot u\,dx, & u\in H^{s_\infty,p}_0(\Omega;\R^n)\\
    \infty,&\text{otherwise}.
\end{cases}$$ Consider the following assumptions over the integrand $W$:
\begin{itemize}
    \item There exists $a\in L^1(\R^n)$ such that $$c_1|A|^p-a(x)\leq W(x,A)\leq a(x)+c_2|A|^p,\,x\in \R^n,\,A\in \R^{n\times n},\,a\in L^1(\R^n).$$
    \item $W(x,\cdot)$ is quasiconvex for a.e. $x\in \Omega$.
\end{itemize} Then we have the following result \cite[Theorem 4.10]{CuetoHidde2025}: 
\begin{Teor}
    The family $(I_\delta)_\delta$, $\Gamma$-converges with respect to the $L^p(\R^n;\R^n)$-convergence to the functional $I_\infty$ as $\delta\to +\infty$,  $$\Gamma(L^p)-\lim_{\delta\to +\infty}I_\delta=I_\infty.$$
\end{Teor}

A direct consequence of this result is the $\Gamma$-convergence of the nonlocal problem with finite horizon \eqref{eqn:Pdelta} to the  problem of fractional linearized elasticity studied in \cite{BellidoCuetoMora2023}, i.e. problem eqref{eqn:P} for $\rho=\rho^{s_\infty}$. In fact, we have that $$\Gamma(L^2)-\lim_{\delta\to +\infty}\begin{cases}
    -\operatorname{div}_{\rho_\delta} \left(C[D_{\rho_\delta,sym}v]\right)&=f,\,\text{in $\Omega$},\\
    \hfill v&=0,\,\text{on $\Omega_{2\delta}\setminus\Omega$}.
\end{cases}=\begin{cases}
    -\operatorname{div}_{s_\infty}\left(C[D^{s_\infty}_{sym}v]\right)&=f,\,\text{in $\Omega$},\\
    \hfill v&=0,\,\text{on $\Omega^c$}.
\end{cases}.$$ 
The equi-coercivity of the family $(I_\delta)_\delta$ and the existence and uniqueness of weak solutions for the nonlocal problem in $H^{\rho,2,\delta}_0(\Omega;\R^n)$ implies the existence and uniqueness of weak solutions for the fractional one in $H^{s_\infty,2}_0(\Omega;\R^n)$.

\subsubsection{Localization in the fractional parameter}
For the particular choice of $D_{\rho^s_\delta}=D_\delta^s$, localization results for the fractional parameter have also been established in\cite{Cueto2023}. For a fixed $\delta>0$,  we define the functional
$$I_s(u):=\begin{cases}
    \displaystyle\int_\Omega W(x,D^s_\delta u)\,dx-\int_\Omega f\cdot u\,dx,&u\in H^{s,p,\delta}_0(\Omega;\R^n),\\
    \infty,&\text{otherwise}.
\end{cases}$$ 
with $f\in L^{p'}(\Omega_\delta;\R^n)$.
 For \(s \to 1^-\), we expect to recover the classical gradient \(D\), since as nonlocality vanishes, the finite horizon no longer influences the behavior of a local operator. Hence, we would expect that the target space of solutions becomes a classical Sobolev space. We define the space $$H^{1,p,\delta}_0(\Omega;\R^n):=\{u\in W^{1,p}_0(\Omega;\R^n): u=0\,\text{a.e. in}\,\Gamma_{\pm\delta}\},$$ and therefore we define the functional $$I_1(u):=\begin{cases}
    \displaystyle\int_\Omega W(x,D u)\,dx-\int_\Omega f\cdot u\,dx,\,&u\in H^{1,p,\delta}_0(\Omega;\R^n)\\
    \infty,\,&\text{else}.
\end{cases}$$


In \cite[Theorem 7]{Cueto2023} it was proved the following $\Gamma$-convergence result.
\begin{Teor}
    Let the family of functionals $(I_s)_s$ and $I_1$ as above, and suppose that $W:\Omega\times\R^{n\times n}\to \R$ is a Carath\'eodory integrand such that there exists two positive constants $c_1,c_2$ with $$c_2|A|^p-c_1\leq W(x,A)\leq C_1(1+|A|^p),$$ for a.e. $x\in \Omega$ and $A\in R^{n\times n}$. Furthermore, we assume that $W(x,\cdot)$ is quasiconvex a.e. in $x\in \Omega_{-\delta}$. Then, the family $(I_s)_s$ $\Gamma$-converges to $I_1$ as $s\to 1^-$ in the sense of $L^p(\Omega_\delta;R^n)$. Moreover, the family $(I_s)_s$ is equi-coercive with respect to convergence in $L^p(\Omega_\delta;\R^n)$.
\end{Teor}

A straightforward consequence of the previous result is 
$$\Gamma(L^2)-\lim_{s\to 1^-}\begin{cases}
    -\operatorname{div}^s_{\delta}(C[D^s_{\delta,sym}v])&=f,\,x\in \Omega_{-\delta},\\
    v&=0,\,x\in \Gamma_{\pm\delta}
\end{cases}=\begin{cases}
    -\operatorname{div}(C[D_{sym}v])&=f,\,x\in \Omega_{-\delta},\\
    v&=0,\,x\in \partial \Omega_{-\delta}
\end{cases},$$ and the existence of weak solution for the minimizers problem  and equi coervicity of the family $I_s$ implies the existence of weak solution for the localized problem.

\subsection{Localization for the Neumann problem}

We also consider localization for the Neumann problem. Again, we follow the ideas in \cite{Cueto2023,CuetoHidde2025, Kreisbeck2024}, but contrary to the previous section, we developed in this section a case which is not considered in those references, localization as $\delta\to 0^+$ for the Neumann problem (the case $s\to 1^-$ is addressed in \cite{Kreisbeck2024}). More concretely, we study the limit problem of the sequence
$$I_\delta(u):=\begin{cases}
    \displaystyle\int_{\Omega}W(x,D_\delta ^su)\,dx-\int_{\Omega_\delta}f\cdot u\,dx, & u\in N^{s,p,\delta}(\Omega;\R^m)^\perp,\\
    \infty, & \text{otherwise},
\end{cases},$$ as $\delta\to 0^+$.

A first key mathematical fact to be taken into account is that for a given domain $\Omega\subset \mathbb{R}^n$ and every $u\in W^{1,p}$, the following convergence 
$$\chi_{\Omega_{\delta}} D_\delta^s u\to Du,$$
as $\delta \to 0^+$ (\cite[Lemma 3.1]{CuetoHidde2025}). Another important mathematical fact concerns the limit of linear manifolds $N^{s,p,\delta}(\Omega)^\perp$ as $\delta \to 0^+$. Since the spaces $N^{s,p,\delta}(\Omega)^\perp$ are the spaces of functions $u\in H^{s,p,\delta}(\Omega)$ such that $\pi_\delta^s(u)=0$, we would expect that the space where the limit problem is formulated is something like $$N^{s,p,0}(\Omega)^\perp:=\{u\in H^{s,p,0}(\Omega):\pi_0^s(u)=0\}.$$ Since $D_0^s=D$, we can identify $H^{s,p,0}(\Omega)$ with $W^{1,p}(\Omega)$, the classical Sobolev space, and $\pi_0^s$ would be the projection from $L^p(\Omega)$ into $N^{s,p,0}(\Omega)$, which again attending at the definition of $N^{s,p,\delta}$ spaces, is just the space of functions in $W^{1,p}(\Omega)$ such that $Du=0$ a.e. in $\Omega$, i.e., functions constant a.e. in $\Omega$. Hence, $\pi_0^s(u)=0$ implies that $$\operatorname{arg}\operatorname{min}_{a\in \R}\norm{u-a}_{L^p(\Omega)}=0,$$ which is equivalent to having $$\int_\Omega |u|^{p-1}\operatorname{sgn}(u)\,dx=0,$$ and hence $$N^{s,p,0}(\Omega)^\perp=\bigg\{u\in W^{1,p}(\Omega): \int_\Omega |u|^{p-1}\operatorname{sgn}(u)\,dx=0\bigg\}.$$ Observe that when $p=2,$ the integral condition is just the usual zero mean value of $u$ in $\Omega$ typically considered for the classical Neumann problem.

Now, we need to establish the Poincar\'e-Wirtinguer inequality of Lemma \ref{PW} with a constant independent of the parameter $\delta$. This is entirely based on the fact that the operator $\mathcal{P}_\delta^s$ has an operator norm independent of $\delta$.\\

Let $\overline{\delta}>0$ and $(\delta_1,\delta_2)\in [0,1]\times [\overline{\delta},1]$. We define the operators $$m_{\delta_1,\delta_2}(\xi):=\frac{\widehat{Q}^s_{\delta_2}(\xi)}{\widehat{Q}_{\delta_1}^s(\xi)},\,M_{\delta_1,\delta_2}(\xi):=\frac{\widehat{\mathcal{P}}^s_{\delta_1}(\xi)}{\widehat{\mathcal{P}}_{\delta_2}^s(\xi)},\,\xi\in \R^n.$$ Observe that by definition of the operator $\mathcal{P}_\delta^s$, $$M_{\delta_1,\delta_2}(\xi)=m_{\delta_1,\delta_2}(\xi),$$ and hence an analogous argument as \cite[Lemma 8]{Cueto2023} (see also \cite[Theorem 3.3]{CuetoHidde2025}) shows that $M_{\delta_1,\delta_2}$ is an $L^p$-multiplier with constant depending on $n,p,\overline{\delta}$. This means that the operator norm of $\mathcal{P}_\delta^s$ is uniformly bounded for the parameter $\delta>0$ and hence we the following lemma, analogous to \cite[Lemma 4.7]{Kreisbeck2024}.

\begin{Lema}\label{PW2}
    Let $1<p<\infty$. Then, there exists a constant $C=C(n,p,\Omega)>0$ such that $$\norm{u-\pi_\delta^s(u)}_{L^p(\Omega_\delta)}\leq C\norm{D_\delta^s u}_{L^p(\Omega;\R^n)},$$ for all $u\in H^{s,p,\delta}(\Omega)$. In particular, if $u\in N^{s,p,\delta}(\Omega)^\perp,$ $$\norm{u}_{L^p(\Omega_\delta)}\leq C\norm{D_\delta^s u}_{L^p(\Omega;\R^n)}.$$
\end{Lema}

As a consequence of the previous lemma we obtain the following result. 
\begin{Lema}\label{Gamma}
    Let $1<p<\infty$ and $(\delta_j)_j\subset (0,1]$ a sequence of horizon parameters such that $\Omega_{-\delta_j}$ is well defined for every $j\in \mathbb{N}$ and such that $\delta_j\to 0$ as $j\to\infty$. Then, it holds that:
    \begin{itemize}
        \item[i)] For all $v\in W^{1,p}(\R^n)$, $\mathcal{P}_{\delta_j}^sv\to v$ in $L^p(\Omega)$ as $j\to \infty$.
        \item[ii)] Let $(u_j)_j$ a sequence of functions such that $u_j\in N^{s,p,\delta_j}(\Omega)^\perp$ for every $j$, satisfying $$\operatorname{sup}_j\norm{D^s_{\delta_j}u_j}_{L^p(\Omega;\R^n)}<\infty.$$ Then, there exists a non relabeled subsequence and a function $u\in N^{s,p,0}(\Omega)^\perp$ such that $u_j\to u$ in $L^p(\Omega)$ and $D^s_{\delta_j}u_j\rightharpoonup Du$ in $L^p(\Omega;\R^n)$, as $j\to \infty$. 
    \end{itemize}
\end{Lema}
\noindent{}\textbf{Proof:}
\begin{itemize}
    \item[i)] Recall that for every $j$, $\mathcal{P}_{\delta_j}^s:W^{1,p}(\R^n)\to H^{s,p,\delta_j}(\R^n)$ is bounded with norm independent of $\delta_j$. Also, the spaces $H^{s,p,\delta_j}(\R^n)$ and $H^{s,p}(\R^n)$ are the same by \cite[Lemma 5]{Cueto2023} (see also (3.4) in \cite{CuetoHidde2025}). Now, let $\overline{\delta}\geq \operatorname{sup}_j\{\delta_j\}_{j\in\mathbb{N}}$. We have that $$\operatorname{sup}_j\norm{\mathcal{P}_{\delta_j}^sv}_{H^{s,p,\overline{\delta}}(\R^n)}\approx \operatorname{sup}_j\norm{\mathcal{P}_{\delta_j}^sv}_{H^{s,p,\delta_j}(\R^n)}<\infty,$$ and since we have the compact embedding $H^{s,p,\overline{\delta}}(\R^n)=H^{s,p}(\R^n)\xhookrightarrow{}\xhookrightarrow{}L^p(\Omega)$, see \cite[Theorem 1.1]{BellidoGarcia2025}, there exists a not relabeled subsequence $\{\mathcal{P}_{\delta_j}^sv\}_{j}$ converging to some $w\in L^p(\Omega)$ in the $p$-norm.\\
    
    We have to prove that $w=v$ in $L^p(\Omega)$. Let $\varphi\in C_c^\infty(\Omega)$. Since $\mathcal{Q}_{\delta_j}^s\varphi\to \varphi$ uniformly as $j\to\infty$ as it was shown in \cite[Lemma 3.1 (i)]{CuetoHidde2025}, we have that $$\int_\Omega w\varphi\,dx=\lim_{j\to\infty}\int_{\Omega_{\delta_j}} (\mathcal{P}_{\delta_j}^sv)(\mathcal{Q}_{\delta_j}^s\varphi)\,dx=\lim_{j\to\infty}\int_{\Omega_{\delta_j}}\left(\mathcal{Q}_{\delta_j}^s\mathcal{P}_{\delta_j}^sv\right)\varphi\,dx=\int_{\Omega}v\varphi\,dx,$$ so $v=w$ in $\Omega$.
    
    \item[ii)] By \eqref{PW2}, $u_j$ is bounded in $L^p(\Omega_{\delta_j})$ for every $j$. The extension operator $\mathcal{E}_{\delta_j}^s:=\mathcal{P}_{\delta_j}^s\circ\mathcal{E}\circ\mathcal{Q}_{\delta_j}^s$, where $\mathcal{E}$ is some extension operator for $W^{1,p}(\Omega)$, is uniformly bounded on $\delta_j$, hence $$\sup_j\norm{\mathcal{E}_{\overline{\delta}}^su_j}_{H^{s,p}(\R^n)}\approx \sup_j\norm{\mathcal{E}_{\overline{\delta}}^su_j}_{H^{s,p,\delta_j}(\R^n)}<\infty,$$ and by the compactness of the embedding of Bessel potential spaces into $L^p(\Omega)$, we can extract a non relabeled subsequence and $w\in L^p(\Omega)$ such that $\mathcal{E}_{\delta_j}^su_j\to w$ in $L^p(\Omega)$, and $D_{\delta_j}^su_j=D_{\delta_j}^s\mathcal{E}_{\delta_j}^su_j\to V$ in $L^p(\Omega;\R^n)$. Observe that given $\psi\in C_c^\infty(\Omega;\R^n)$, we have that $$\int_\Omega V\cdot \psi\,dx=\lim_{j\to\infty}\int_{\Omega}D_{\delta_j}^su_j\cdot \psi\,dx=-\lim_{j\to\infty}\int_{\Omega_{\delta_j}}u_j\operatorname{div}_{\delta_j}^s\psi\,dx=-\int_\Omega u\operatorname{div}\psi\,dx,$$ thus $w\in W^{1,p}(\Omega)$, and $D^s_{\delta_j}u_j\rightharpoonup Du$ in $L^p(\Omega;\R^n)$, as $j\to\infty$. Since $u_j-\mathcal{E}_{\delta_j}^su_j\in N^{s,p,\delta_j}(\Omega)$, we can write $$u_j=\mathcal{E}_{\delta_j}^su_j+\left(u_j-\mathcal{E}_{\delta_j}^su_j\right)=\mathcal{E}_{\delta_j}^su_j-\pi_{\delta_j}^s\left(\mathcal{E}_{\delta_j}^su_j\right).$$ Adapting the proof of \cite[Lemma 6.3 (ii)]{Kreisbeck2024}, we obtain that $\pi_{\delta_j}^s\left(\mathcal{E}_{\delta_j}^su_j\right)\to \pi_0^s(w)$ in $L^p(\Omega)$ as $j\to \infty$. Hence, $$u_j=\mathcal{E}_{\delta_j}^su_j-\pi_{\delta_j}^s\left(\mathcal{E}_{\delta_j}^su_j\right)\to w-\pi_0^s(w),$$ as $j\to \infty$. If we define $u:=w-\pi_0^s(w)$, $u_j\to u$ as $j\to \infty$, and $u\in N^{s,p,0}(\Omega)^\perp$, since $\pi_0^s(w)$ is constant in $\Omega$ due to the fact that $D_{\delta_j}^su_j\rightharpoonup Dw$, in $L^p(\Omega)$, i.e., $D\pi_0^s(w)=0$ in $\Omega$.
\end{itemize}\qed \\

Now, we are ready to establish the $\Gamma$-convergence result for the case of vanishing horizon, $\delta\to 0^+$.
\begin{Teor}
    Let $1<p<\infty$, $f\in L^{p'}(\Omega_\delta;\R^m)$ and $W:\Omega\times \R^{m\times n}\to \R\cup\{\infty\}$ be a Carath\'eodory integrand such that $$W(x,A)\geq c\left(|A|^p-1\right),$$ with a constant $c>0$, for a.e. $x\in \Omega$ and for every $A\in \R^{m\times n}$. If $v\mapsto \int_\Omega W(x,Dv)\,dx$ is weakly lower semicontinuous on $W^{1,p}(\Omega;\R^m)$, then the family of functionals $(I_\delta)_\delta$ with $I_\delta:L^p(\Omega_\delta;\R^m)\to \R\cup \{\infty\}$, defined by $$I_\delta(u):=\begin{cases}
   \displaystyle \int_{\Omega}W(x,D_\delta ^su)\,dx-\int_{\Omega_\delta}f\cdot u\,dx, & u\in N^{s,p,\delta}(\Omega;\R^n)^\perp,\\
    \infty, & \text{otherwise},
\end{cases}$$ $\Gamma$-converges with respect to $L^p(\Omega;\R^m)$-convergence as $\delta\to 0^+$ to $I_0:L^p(\Omega;\R^n)\to \R\cup\{\infty\}$, defined as $$I_0(u):=\begin{cases}
    \displaystyle\int_{\Omega}W(x,Du)\,dx-\int_{\Omega}f\cdot u\,dx, & u\in N^{s,p,0}(\Omega;\R^n)^\perp,\\
    \infty, &\text{otherwise}
\end{cases}$$
Also, the family $(I_\delta)_\delta$ is equi-coercive in $L^p(\Omega;\R^m)$.
\end{Teor}
\noindent{}\textbf{Proof:} The first step is to prove the equi-coercivity of the family of functionals, which follows directly from the lower bound $f(x,A)\geq c(|A|^p-1),$ the Nonlocal Poincar\'e-Wirtinguer inequality with constant independent of the horizon parameter of Lemma \ref{PW2} and Lemma \ref{Gamma} ii). All this together implies that for every sequence of functions $(u_j)_j$ such that $u_j\in N^{s,p,\delta_j}(\Omega;\R^m)^\perp$ with $\operatorname{sup}_jI_{\delta_j}(u_j)<\infty$, admits a subsequence converging strongly in the $L^p(\Omega;\R^m)$ to a function $u\in N^{s,p,0}(\Omega;\R^m)^\perp.$\\

The second step is to prove the \textit{Liminf-inequality}. Let $(u_j)_j$ a family of functions such that $u_j\in N^{s,p,\delta_j}(\Omega;\R^m)^\perp$ for every $j$ such that $u_j\to u$ in $L^p(\Omega;\R^m)$, and we assume without loss of generality that $\operatorname{sup}_jI_{\delta_j}(u_j)<\infty$. The lower bound $$W(x,A)\geq c(|A|^p-1),$$ 
together with Lemmas \ref{PoincareII} and  \ref{Gamma} {\it (ii)} implies that $u\in N^{s,p,0}(\Omega;\R^m)^\perp$ and $D_{\delta_j}^su_j\rightharpoonup Du$ in $L^p(\Omega;\R^{m\times n})$ as $j\to \infty$. Now, we define the sequence $v_j:=\mathcal{Q}_{\delta_j}^s u_j$ for every $j$. Note that $v_j\in W^{1,p}(\Omega)$ with $v_j\to u$ and $Dv_j=D_{\delta_j}^su_j\to Du$ as $j\to \infty$. The result follow by the the weak lower semicontinuity of $v\mapsto \int_\Omega W(x,Dv)\,dx$, we got that \begin{align*}\liminf_{j\to\infty}I_{\delta_j}(u)&=\liminf_{j\to\infty}\left(\int_\Omega W(x,D_jv)\,dx-\int_{\Omega_{\delta_j}}f\cdot u_j\,dx\right)\\
&\geq \int_\Omega W(x,Du)\,dx-\int_\Omega f\cdot u\,dx=I_0(u).\end{align*}
The last step is the existence of recovery sequence. Let $u\in N^{s,p,0}(\Omega;\R^m)^\perp$ with $I_0(u)<\infty$. Take $v\in W^{1,p}(\R^n;\R^m)$ such that $v|_\Omega=u$. Exploiting the fact that $\mathcal{P}_{\delta_j}^s:W^{1,p}(\R^n;\R^m)\to H^{s,p,\delta_j}(\R^n;\R^m)$ continuously, we define the sequence $$u_j:=\left(\mathcal{P}_{\delta_j}^sv\right)\bigg|_{\Omega_{\delta_j}}-\pi_{\delta_j}^s\left(\left(\mathcal{P}_{\delta_j}^s v\right)\bigg|_{\Omega_{\delta_j}}\right)\in N^{s,p,\delta_j}(\Omega;\R^m)^\perp.$$ By construction, $$D_{\delta_j}^su_j=D\mathcal{Q}_{\delta_j}^s\left(\mathcal{P}_{\delta_j}^sv\right)\bigg|_{\Omega_{\delta_j}}=Dv|_{\Omega}=Du.$$ Also, by Lemma \ref{Gamma}, $u_j\to v-\pi_0^s(v)$ in $L^p(\Omega;\R^m)$ as $j\to \infty$. Note that $v\in N^{s,p,0}(\Omega;\R^m)^\perp$ and $Dv=Du$, hence we can identify $u$ with $v-\pi_0^s(v)$ in $\Omega$. Finally, we observe that \begin{align*}
    I_{\delta_j}(u_j)&=\int_{\Omega}W(x,D_{\delta_j}^su_j)\,dx-\int_{\Omega_{\delta_j}}f\cdot u_j\,dx=\int_\Omega W(x,Du)\,dx-\int_{\Omega_{\delta_j}}f\cdot u_j\,dx\to I_0(u),
\end{align*} as $j\to \infty$. \qed \\

Now that we have established the $\Gamma$-convergence of the problem $I_\delta$, we directly obtain the convergence of the problem \eqref{eqn:PN} into the linear elasticity problem with Neumann boundary conditions. In particular, \begin{equation*}\Gamma(L^p)-\lim_{\delta\to 0^+}\begin{cases}
    -\operatorname{div}_\delta^s(C[D_{\delta,sym}^su])=f,\,&\text{in $\Omega_{-\delta}$},\\
    \mathcal{N}_\delta^s(CD_{\delta,sym}^su)=0,\,&\text{on $\Gamma_{\pm\delta}$}
\end{cases}=\begin{cases}
    -\operatorname{div}(C[D_{sym}u])=f,\,&\text{in $\Omega$},\\
    (C[D_{sym}u])\cdot \nu=0,\,&\text{on $\partial\Omega$}
\end{cases},\end{equation*} where $\nu$ is an outward pointing unit normal to $\partial\Omega$. Moreover, due to the equicoercivity of the family $I_\delta$, the minimizers of the problem $I_\delta$ converge up to a subsequence in $L^p(\Omega;\R^n)$ to a weak solution for the local elasticity problem when $\delta\to 0^+$.

\section*{Acknowledgements}
This work was supported by {\it Agencia Estatal de Investigación} (Spain) through grant PID2023-151823NB-I00 and {\it Junta de Comunidades de Castilla-La Mancha} (Spain) through grant SBPLY/23/180225/000023. G.G.-S. is supported by a Doctoral Fellowship by Universidad de Castilla-La Mancha \text{2024-UNIVERS-12844-404}.

\section*{Conflicts of interest}
The authors declare that there are no conflicts of interest regarding the publication of this paper.

\addcontentsline{toc}{section}{References}
\bibliographystyle{plain} 

\end{document}